\documentclass[aos]{imsart}

\RequirePackage{amsthm,amsmath,amsfonts,amssymb}
\RequirePackage[numbers,sort&compress]{natbib}
\RequirePackage[colorlinks,citecolor=blue,urlcolor=blue]{hyperref}
\RequirePackage{graphicx}
\RequirePackage{enumitem}
\RequirePackage{mathrsfs}

\startlocaldefs
\theoremstyle{plain}

\newtheorem{theorem}{Theorem}[section]
\newtheorem{lemma}[theorem]{Lemma}
\theoremstyle{remark}
\newtheorem{remark}{Remark}[section]

\newcommand{\N}{\mathbb{N}}

\newcommand{\R}{\mathbb{R}}
\newcommand{\I}{\mathbb{I}}

\newcommand{\SgE}{\mathcal{E}}
\renewcommand{\P}{\mathbb{P}}
\newcommand{\E}{\mathbb{E}}

\newcommand{\atom}{\boldsymbol{\alpha}}
\newcommand{\chain}{\textbf{X}}
\newcommand{\splitchain}{{\check{\chain}}}
\newcommand{\measures}[1][\SgE]{\mathscr{M}{\left( {#1} \right)_+}}
\newcommand{\measurables}[1][\pi]{L^1\left(E,{#1}\right)}

\newcommand{\block}{\mathcal{B}}

\newcommand{\numreg}[1][n]{T\left({#1}\right)}
\newcommand{\irreducibilityMeasure}{\psi}
\renewcommand{\leq}{\leqslant}
\renewcommand{\geq}{\geqslant}

\newcommand{\expectation}[2]{\E_{#1}{\left(#2\right)}}

\def\argmax{\mathop{\mbox{\sl\em argmax}}}

\renewcommand{\epsilon}{\varepsilon}
\renewcommand{\hat}{\widehat}
\renewcommand{\tilde}{\widetilde}

\newcommand{\Tn}{T_n(C)}

\newcommand{\asymptnumbreg}{u\left(n\right)}
\newcommand{\distfunc}{F}

\newcommand{\AssumpXText}{(A1)}
\newcommand{\AssumpX}{\hyperref[assumpt: assumpX]{\AssumpXText}}
\newcommand{\AssumpWText}{(A2)}
\newcommand{\AssumpW}{\hyperref[assumpt: assumpW]{\AssumpWText}}
\newcommand{\AssumpCText}{(A3)}
\newcommand{\AssumpC}{\hyperref[assumpt: assumpC]{\AssumpCText}}
\newcommand{\AssumpInteriorPointText}{(A4)}
\newcommand{\AssumpInteriorPoint}{\hyperref[assumpt: assumpInteriorPoint]{\AssumpInteriorPointText}}
\newcommand{\FincreasingText}{(A5)}
\newcommand{\Fincreasing}{\hyperref[assumpt: Fincreasing]{\FincreasingText}}
\newcommand{\fdecreasingText}{(A6)}
\newcommand{\fdecreasing}{\hyperref[assumpt: fdecreasing]{\fdecreasingText}}
\newcommand{\fcontinuousText}{(A7)}
\newcommand{\fcontinuous}{\hyperref[assumpt: fcontinuous]{\fcontinuousText}}
\newcommand{\AssumpXrateText}{(B1)}
\newcommand{\AssumpXrate}{\hyperref[assumpt: AssumpXrate]{\AssumpXrateText}}
\newcommand{\AssumpDeriveeText}{(B2)}
\newcommand{\AssumpDerivee}{\hyperref[assumpt: AssumpDerivee]{\AssumpDeriveeText}}
\newcommand{\AssumpInitialMeasureText}{(B3)}
\newcommand{\AssumpInitialMeasure}{\hyperref[assumpt: AssumpInitialMeasure]{\AssumpInitialMeasureText}}
\newcommand{\AssumpSizeBlockText}{(B4)}
\newcommand{\AssumpSizeBlock}{\hyperref[assumpt: AssumpSizeBlock]{\AssumpSizeBlockText}}

\endlocaldefs

\begin{document}

\begin{frontmatter}
  \title{Harris recurrent Markov chains and nonlinear monotone cointegrated models}
  \runtitle{Harris Markov chains and nonlinear monotone cointegrated models}

  \begin{aug}
    \author[A]{\fnms{Patrice}~\snm{Bertail}\ead[label=e1]{patrice.bertail@parisnanterre.fr}\orcid{0000-0002-6011-3432}}
    \author[A]{\fnms{Cécile}~\snm{Durot}\ead[label=e2]{cecile.durot@parisnanterre.fr}\orcid{0000-0002-4502-0174}}
    \and
    \author[A,B]{\fnms{Carlos} \snm{Fernández}\ead[label=e3]{fernandez@telecom-paris.fr}\orcid{0000-0002-8577-865X}}
    \address[A]{MODAL'X, UMR CNRS 9023, Universit\'e Paris Nanterre\printead[presep={,\ }]{e1,e2}}

    \address[B]{LTCI,
      Telecom Paris, Institut Polytechnique de Paris\printead[presep={,\ }]{e3}}
  \end{aug}

  \begin{abstract}
    In this paper, we study a nonlinear cointegration-type model of the form
    \(Z_t = f_0(X_t) + W_t\) where \(f_0\) is a monotone function and \(X_t\)
    is a Harris recurrent Markov chain. We use a nonparametric Least Square
    Estimator to locally estimate \(f_0\), and under mild conditions, we show
    its strong consistency and obtain its rate of convergence. New results
    (of the Glivenko-Cantelli type) for localized null recurrent Markov chains
    are also proved.
  \end{abstract}



\end{frontmatter}

\section{Introduction and motivations}\label{sec: introduction}
\subsection{Linear and nonlinear cointegation models}
\textit{Linear cointegration} introduced by  \cite{Granger1981}  and developed by \cite{EngleGranger1987} and \cite{Johansen1988,Johansen1991}, is a concept used in statistics and econometrics to describe a long-term relationship between two or more time series. In general, these time series are non-stationary, integrated of order 1, that is, they behave roughly as random walks. In traditional linear cointegration analysis, variables are assumed to have a linear relationship, which means their long-term equilibrium, as time grows, is characterized by a constant linear combination. This concept has since been extensively studied, particularly in the field of econometrics \cite{EngleGranger1987, Phillips1991,PhillipsSolo1992, Johansen1988,Johansen1991}. Notice that, when there is indeed a significant linear relationship, the link is monotone in each of the variables.

However, in some cases, the relationship between variables may exhibit nonlinear
behavior, which cannot be adequately captured by linear cointegration models. The
incorporation of nonlinearities allows for a more comprehensive understanding of
long-term relationships between variables. \cite{HANSEO2002} have developed an approach
for analyzing nonlinear cointegration through threshold cointegration models. These
models assume that the linear relationship between variables differs after some
changepoints, leading to different long-run equilibrium states (for instance according to
some latent regimes). Threshold cointegration models provide a framework for capturing
nonlinearities in the data and estimating the changepoints. Refer to \cite{stig2020} for
examples and discussions on the importance to introduce nonlinearities in cointegration
applications and for further references.

Another method for analyzing nonlinear cointegration is through the use of smooth
transition cointegration models introduced by \cite{GranTer1993}. These models assume a
ECM (Error Correction Model) form and allow for smooth transitions between different
regimes in the data. Most estimators of non-linear cointegration may be seen as
Nadaraya-Watson estimators of the link function. For instance, Wang and Phillips
\cite{WangPhillips2009a,WangPhillips2009b, Wang2015} show that it is possible to estimate
and perform asymptotic inference in specific nonparametric cointegration regression
models using kernel regression techniques. Furthermore, they established that the
self-normalized kernel regression estimators converge to a standard normal distribution
limit, even when the explanatory variable is integrated. These findings indicate that the
estimators can consistently capture the underlying relationship between variables, even
in cases where the explanatory variable exhibits non-stationary behavior. The problem of
estimating \(f_0\) under the Markovian assumptions has also been tackled using local
linear M-type estimators in \cite{CaisElias2003, ZhengyanLiChen2009} using smoothing
techniques.

These results have been partially extended in the framework of general $\beta$-recurrent
Markov chain and not just integrated $I(1)$ time series by
\cite{karlsen2007nonparametric, CaiTjostheim2015}. Consider the simple framework where we
observe two Markov chains, $Z_t$ and $X_t$. \cite{karlsen2007nonparametric} are
essentially interested in the study of \textit{nonlinear cointegrated} models such as,
\begin{equation}\label{eq: model}
  Z_t = f_0(X_t) + W_t,
\end{equation}
where $f_0$ is a nonlinear function.  $X_t$  and
\(W_t\) are independent processes, and \(X_t\) is a positive or \(\beta\)-null recurrent
Markov chain. Despite the fact that there is no stationary probability measure for $X_t$, they apply the Nadaraya-Watson method to estimate \(f_0\) and
established the asymptotic theory of the proposed estimator. The rate of convergence of the estimators at some point $x$ is essentially linked to the local properties of the $\beta$ -null recurrent chain $X_t$ and typically of the order of the square root of the number of visits of the chain in a neighborhood of the point $x$.

\subsection{Monotone cointegration models: motivations}\label{sec: motivations}

Monotonicity in cointegration is a rather natural assumption in many economic
applications, for instance for modeling demand as a function of income or prices (see for
instance \cite{deamuel1980}) or other variables. Suppose, for example, we are interested
in analyzing the long-term relationship between ice cream sales and the average monthly
temperature. These two non-stationary variables may be modeled by some $\beta$-recurrent
Markov chains. We hypothesize that as the average monthly temperature increases, the
demand for ice cream also increases: however the rate of increase may vary according to
the season. In that case, the nonlinear relationship between the two variables will be
monotone. In microeconomics, the same phenomenon is expected for Engel curves, describing
how real expenditure varies with household income (see \cite{deamuel1980}). Expenditures
and income (or their log in most models) may be considered as non-stationary variables.
However considering a linear cointegration between them may be misleading, since the
relationship may change along the life cycle. By Engle's law, the relationship between
the two variables should be monotone. Other types of examples of monotone non-linear
cointegration phenomenon may be found in \cite{stig2020}.

The purpose of this paper is to propose a simple estimator that is automatically
monotone, does not require strong smoothness assumptions (we only require continuity of
the link function), and operates under general Markovian assumptions. We establish a
nonparametric estimation theory of the nonparametric least squares estimator (LSE) for
the function $f_0$ in the model \eqref{eq: model} under the constraints that $f_0$ is
monotone non-increasing. Here, $\{W_t\}$ is an unobserved process such that
$E(W_t|X_t)=0$ to ensure identifiability of $f_0$. Since a minimal condition for
undertaking asymptotic analysis on \(f_0(x_0)\) at a given point $x_0$ is that, as the
number of observations on $\{X_t\}$ increases, there must be infinitely many observations
in the neighborhood of $x_0$, the process $\{X_t\}$ will be assumed to be a Harris
recurrent Markov chain (cf section \ref{sec: markov_theory}). We consider at the same
time the stationary and $\beta$ null recurrent non-stationary framework. To our
knowledge, it is the first time that such an estimator is proposed in the literature in
such a large framework.

\subsection{The estimator}\label{sec: estimator}

Let \(C\) be a set whose interior contains our point of interest \(x_0\). Having observed
\(\left(X_t,Z_t\right)\}_{t=0}^{n}\), we denote by \(\Tn\) the number of times that
\(\chain\) visited \(C\) up to time \(n\) and by \(\sigma_C\left(i\right)\) the time of
the \(i\)-th visit. Then, we consider the nonparametric LSE defined as the minimizer of

\begin{equation}\label{eq: LSE}
  f \mapsto \sum\limits_{i = 1}^{\Tn} {{{\left( {{Z_{{\sigma_C\left(i\right)}}} - f\left( {{X_{{\sigma_C\left(i\right)}}}} \right)} \right)}^2}}
\end{equation}
over the set of non-increasing functions $f$ on $\R$. The nonparametric LSE
$\hat f_n$  has a well known characterization, as follows. Let $m$ be the number
of unique values of $X_{\sigma_C\left(1\right)},\dots,X_{\sigma_C\left(\Tn\right)}$, and  $Y_{1}<\dots
  <Y_{m}$ be the corresponding order statistics.
Then, $\hat f_n(Y_k)$ is the left-hand slope at
$\sum_{i=1}^{\Tn} \I{\left\{X_{\sigma_C\left(i\right)}\leq Y_k\right\}}$ of the least concave majorant
of the set of points
\begin{equation}\label{eq: CSD}
  \left\{(0,0),\ \left(\sum_{i=1}^{\Tn} \I{\left\{X_{\sigma_C\left(i\right)}\leq Y_k\right\}},\sum_{i=1}^{\Tn} Z_{\sigma_C\left(i\right)}\I{\left\{X_{\sigma_C\left(i\right)}\leq Y_k\right\}}\right), \ k=1,\dots,m
  \right\},
\end{equation}
and it can be computed using simple algorithms as discussed in \citep{barlow1972statistical}.
Thus, the constrained LSE is uniquely defined at the observation points, however, it is not
uniquely defined between these points: any monotone interpolation of these values is a
constrained LSE. As is customary, we consider in the sequel the piecewise-constant and
left-continuous LSE that is constant on every interval
$(Y_{k-1},Y_{k}]$, $k=2,\dots,m$ and also on $(-\infty,Y_{1}]$ and on
$[Y_{m},\infty)$.

The use of a localized estimator is due to the fact that we need to control the behavior
of the chain around \(x_0\), and, to do this, we need to estimate the asymptotic
"distribution" of \(\chain\) in a vicinity of \(x_0\). For Harris recurrent Markov
chains, the long-term behavior of the chain is given by its invariant measure (see
Section \ref{sec: markov_theory}). In the positive recurrent case, the invariant measure
is finite and it can be estimated by simply considering the empirical cumulative
distribution function of the \(X_t\). However, in the null recurrent case, the invariant
measure is only \(\sigma\)-finite, hence, we need to localize our analysis in a set big
enough such that the chain visits it infinitely often, but small enough that the
restriction of the invariant measure to it is finite. Moreover, contrary to the bandwidth
in kernel type estimators, $C$ does not depend on $n$, and the rate of convergence of the
estimator does not depend on $C$.

\subsection{Outline}\label{sec: outline}

Since our paper draws quite heavily on the theory of Harris recurrent Markov chains, we
have added a small introduction to the subject as well as the main results that we use
throughout the paper in Section \ref{sec: markov_theory}. In Section \ref{sec:
consistency}, we show that under very general assumptions, our estimator \(\hat f_n\) is
strongly consistent, while its rate of convergence is presented in Section \ref{sec:
rates}. In Section \ref{sec: localized_markov_chains}, we present new results concerning
the localized empirical process of Harris recurrent Markov chains that have emerged
during our investigation and we believe are interesting in their own right. Section
\ref{sec: proofs} contains the proofs of our main results.

\section{Markov chain theory and notation}\label{sec: markov_theory}

In this section, we present the notation and main results related to Markov chains that
are needed to present our main results. For further details, we refer the reader to the
first section of the Supplementary Material \cite{BertailDurotFernandezSupplements2023}
and the books \cite{Nummelin1984,Meyn2009,markovChain2018}.

Consider a time-homogeneous irreducible Markov Chain, denoted as \(\chain = X_0, X_1,
X_2, \ldots\), defined on a probability space \(\left(E, \SgE,\P\right)\), where \(\SgE\)
is countably generated. The irreducibility measure of the chain is represented by
$\irreducibilityMeasure$. The transition kernel of the chain is denoted as
\(P\left(x,A\right)\) and its initial distribution is represented by $\lambda$. If the
initial measure of the chain is specified, we use \(\P_{\lambda}\) (and \(\E_{\lambda}\))
to denote the probability (and the expectation) conditioned on the law of the initial
state \(\mathcal{L}\left(X_0\right)=\lambda\).

For any set \(C\in\SgE\), we will denote by \(\sigma_C\) and \(\tau_C\), respectively,
the times of first visit and first return of the chain to the set \(C\), i.e.
\(\tau_C=\operatorname{inf}\left\{ n\geq 1: X_n\in C \right\}\) and
\(\sigma_C=\operatorname{inf}\left\{ n\geq 0: X_n\in C \right\}\). The subsequent visit
and return times \(\sigma_C, \tau_C\left(k\right)\), \(k\geq 1\) are defined inductively.

Given that our methods will only deal with the values of \(\chain\) in a fixed set \(C\),
if \(A\) is a measurable set, we will write \(\I_C\{X_t\in A\}\) instead of \(\I\{X_t\in
A\cap C\}\) and if \(A=E\), then we will simply write \(\I_C\left(X_t\right)\). We will
use \(T_n\left(C\right)\) to denote the random variable that counts the number of times
the chain has visited the set \(C\) up to time \(n\), that is \(T_n\left( C
\right)=\sum_{t=0}^{n}{\I_C\left(X_t\right)}\). Similarly, we will write
\(T\left(C\right)\) for the total of numbers of visits the chain \(\chain\) to \(C\). The
set \(C\) is called \textit{recurrent} if \(\E_x T\left(C\right)=+\infty\) for all \(x\in
C\) and the chain \(\chain\) is recurrent if every set $A\in\SgE$ such that
$\irreducibilityMeasure\left(A\right)>0$ is recurrent. A recurrent chain is called
\textit{Harris recurrent} if for all \(x\in E\) and all \(A\in\SgE\) with
\(\irreducibilityMeasure(A)>0\) we have \(\P\left( {{X_n} \in A\; \text{infinitely
often}\left| {{X_0} = x} \right.} \right) = 1\).

Denote by \(\SgE^+\) the class of nonnegative measurable functions with positive
\(\irreducibilityMeasure\) support. A function \(s\in\SgE^+\) is called \textit{small} if
there exists an integer \(m_0\geq 1\) and a measure \(\nu\in\measures\) such that
\begin{equation}\label{eq: minorization_general}
  {P ^{m_0}}\left( {x,A} \right) \geqslant s\left( x \right)\nu \left( A \right)\quad\forall x\in E, A\in\SgE.
\end{equation}
When a chain possesses a small function \(s\), we say that it satisfies the
\textit{minorization inequality} \(M\left(m_0,s,\nu\right)\). A set \(A\in\SgE\) is said to be \textit{small} if the function \(\I_A\) is small.
Similarly, a measure \(\nu\) is \textit{small} if there exist \(m_0\), and \(s\) that
satisfy \eqref{eq: minorization_general}. By Theorem 2.1 in \cite{Nummelin1984}, every
irreducible Markov chain possesses a small function and Proposition 2.6 of the same book
shows that every measurable set \(A\) with \(\irreducibilityMeasure\left(A\right)>0\)
contains a small set. In practice, finding such a set consists in most cases in
exhibiting an accessible set, for which the probability that the chain returns to it in
$m$ steps is uniformly bounded below. Moreover, under quite wide conditions a compact set
will be small, see \cite{Feigin1985}.

An irreducible chain possesses an accessible atom, if there is a set \(\atom\in\SgE\)
such that for all \(x,y\) in \(\atom\): \(P(x,.)=P(y,.)\) and
\(\irreducibilityMeasure(\atom)>0\). When an accessible atom exists, the
\textit{stochastic stability} properties of \(\chain\) amount to properties concerning
the speed of return time to the atom only. Moreover, it follows from the \textit{strong
Markov property} that the sample paths may be divided into independent blocks of random
length corresponding to consecutive visits to \(\atom\). The sequence
\(\left\{\tau_{\atom}(j)\right\} _{j\geqslant 1}\) defines successive times at which the
chain forgets its past, called \textit{regeneration times}. Similarly, the sequence of
i.i.d. blocks \(\{\block_j\}_{j\geq 1}\) are named \textit{regeneration blocks}. The
random variable \(\numreg=T_n\left(\atom\right)-1\) counts the number of i.i.d. blocks up
to time \(n\). This term is called \textit{number of regenerations up to time \(n\)}.

If \(\chain\) does not possess an atom but is Harris recurrent (and therefore satisfies a
minorization inequality \(M\left(m_0,s,\nu\right)\)), a \textit{splitting technique},
introduced in \cite{Nummelin1978, Nummelin1984}, allows us to extend in some sense the
probabilistic structure of \(\chain\) in order to artificially construct an atomic chain
(named the \textit{split chain} and denoted by $\splitchain$) that inherits the
communication and stochastic stability properties from $\chain$. One of the main results
derived from this construction is the fact that every Harris recurrent Markov chain
admits a unique (up to multiplicative constant) invariant measure (see Proposition 10.4.2
in \cite{Meyn2009}), that is, a measure $\pi$ such that

\[\pi \left( B \right) = \int {P\left( {x,B} \right)d\pi \left( x \right)}.\quad\forall B\in\SgE.\]

The invariant measure is finite if and only if
\(\E_{\check{\atom}}{\tau_{\check{\atom}}}<+\infty\), in this case we say the chain is
\textit{positive recurrent}, otherwise, we say the chain is \textit{null recurrent}. A
null recurrent chain is called \(\beta\)-null recurrent (c.f. Definition 3.2 in
\cite{Tjostheim-2001}) if there exists a small nonnegative function \(h\), a probability
measure \(\lambda\), a constant \(\beta \in \left( {0,1} \right)\) and a slowly varying
function \(L_h\) such that
\[{\E_\lambda }\left( {\sum\limits_{t = 0}^n {h\left( {{X_t}} \right)} } \right) \sim \frac{1}{{\Gamma \left( {1 + \beta } \right)}}{n^\beta }{L_h}\left( n \right)\quad\textnormal{as}\;n\to\infty.\]

As argued in \cite{Tjostheim-2001}, is not a too severe restriction to assume \(m_0=1\).
Therefore, throughout this paper we assume that \(\chain\) satisfies the minorization
inequality \(M(1,s,\nu)\), i.e, there exist a measurable function \(s\) and a probability
measure \(\nu\) such that \(0\leq s\left(x\right)\leq 1\), \(\int_{E}{s(x)d\nu(x)}>0\)
and
\begin{equation}\label{eq: minorization}
  {P}\left( {x,A} \right) \geqslant s\left( x \right)\nu \left( A \right).
\end{equation}

\begin{remark}
  The extensions to the case where \(m_0>1\) of the results that will be presented in this paper
  can be carried out (although they involve some complicated notations/proofs) using the
  \(m\)-skelethon or the resolvent chains, as described in \cite{Chen1999, Chen2000} and
  Chapter 17 of \cite{Meyn2009}. However, they are not treated in this paper.
\end{remark}

The following theorem is a compendium of the main properties of Harris's recurrent Markov
chains that will be used throughout the paper. Among other things, it shows that the
asymptotic behavior of \(\numreg\) is similar to the function \(u\left(n\right)\) defined
as
\begin{equation}\label{eq: definition_un}
  u\left(n\right)=\begin{cases}
    n,                       & \text{if }\chain \text{ is positive recurrent}          \\
    n^\beta L\left(n\right), & \text{if }\chain \text{ is }\beta\text{-null recurrent}
  \end{cases}.
\end{equation}

\begin{theorem}\label{th: t_n_general_results} Suppose \(\chain\) is a Harris recurrent, irreducible Markov chain, with initial measure \(\lambda\), that
  satisfies the minorization condition \eqref{eq: minorization}. Let \(\numreg\)
  be the number of complete regenerations until time \(n\) of the split chain \(\splitchain\)
  , let \(C\in\SgE\) be a small set and \(\pi\) be an invariant measure for \(\chain\).
  Then,
  \begin{enumerate}
    \item \(0<\pi\left(C\right)<+\infty\).
    \item \(\frac{\numreg}{\Tn}\) converges almost surely to a positive constant.
    \item \(\frac{\numreg}{u\left(n\right)}\) converges almost surely to a positive
          constant if \(\chain\) is positive recurrent
          and converges in distribution to a Mittag-Leffler\footnote{The Mittag-Leffler  distribution with index \(\beta\)
            is a non-negative continuous distribution, whose moments
            are given by
            \begin{equation*}
              \E \left(M^m_\beta\left(1\right)\right)=\frac{m!}{\Gamma\left(1+m\beta\right)}\;\;m\geq 0.
            \end{equation*}
            See (3.39) in \cite{Tjostheim-2001} for more details.}
          random variable with index \(\beta\) if \(\chain\) is \(\beta\)-null recurrent.
  \end{enumerate}
\end{theorem}

\section{Consistency}\label{sec: consistency}

The aim of the section is to show that for an arbitrary $x_0$ in the support of $f_0$,
the LSE $\hat f_n(x_0)$ is consistent. We make the following assumptions on the processes
$\chain=\{X_t\}$ and $\textbf{W}=\{W_t\}$.

\begin{itemize}
  \item[\AssumpXText]\label{assumpt: assumpX}
        $\chain$ is a Harris recurrent Markov chain whose kernel \(P(x,A)\)
        satisfies the minorization condition (\ref{eq: minorization}).
\end{itemize}

Let \(\mathcal{F}_n=\sigma\left(\left\{ X_0, \ldots, X_n \right\}\right)\) be sigma
algebra generated by the chain \( \chain \) up to time \(n\).

\begin{itemize}
  \item[\AssumpWText]\label{assumpt: assumpW} \text For each \(n\), the random variables
        $W_1,\ldots,W_n$ are conditionally independent given ${{\mathcal{F}_n}}$,
        \(\E{\left(W_t|\mathcal{F}_n\right)}=0\) and \(\operatorname{Var}\left(W_t|\mathcal{F}_n\right)\leq\sigma^2\) for some $\sigma>0$.
\end{itemize}

It follows from Assumption {\AssumpX} that the Markov Chain $\chain$ admits a unique (up
to a multiplicative constant) $\sigma$-finite invariant measure $\pi$. Let $C$ be a set
such that $0<\pi\left(C\right)<\infty$ and $x_0\in C$. We denote by $F_n$ the process
defined by
\begin{equation}\label{eq: Fn}
  F_n(y)=\frac 1{\Tn} \sum_{i=1}^{\Tn}{\I{\{ X_{\sigma_C\left(i\right)}\leq y \}}} = \frac 1{\Tn}\sum_{t=0}^n\I_C{\{X_t\leq y\}}
\end{equation}
for all $y\in\R$, which is a localized version of the empirical distribution function of
the $X_t$'s. It is proved
in Lemma \ref{lem: cvFn} that \(F_n\) converges almost surely to the distribution function $\distfunc$ supported on $C$
and defined by
\begin{equation}\label{eq: F}
  {\distfunc}\left( y \right) = \frac{{\pi \left( {C \cap \left( { - \infty ,y} \right]} \right)}}{{\pi \left( C \right)}}.
\end{equation}

Our next two assumptions guarantee that there is a compact \(C\), that is a small set and
contains \(x_0\) as an interior point. Sets like this can be found under very wide
conditions (cf \cite{Feigin1985}).

\begin{itemize}
  \item[\AssumpCText]\label{assumpt: assumpC} There is $\delta=\delta(x_0)$ such that the set
        $C = \left[ {{x_0} - \delta ,{x_0} + \delta } \right]$
        is \textit{small}.
  \item[\AssumpInteriorPointText]\label{assumpt: assumpInteriorPoint} $x_0$ belongs to the interior of the support of \(X_t\).
\end{itemize}

Notice that by part 1 of Theorem \ref{th: t_n_general_results}, \AssumpC$\,$ guarantees
that \(\pi(C)\) is finite and positive, and hence, \(F\) is properly defined.

In addition to the assumptions on the processes $\{X_t\}$ and $\{W_t\}$, we need
smoothness assumptions on $\distfunc$ and on $f_0$. In particular, we will assume that
$\distfunc$ and $f_0$ are continuous and strictly monotone in $C$. This implies that
$f_0$ and $\distfunc$ are invertible in $C$, so we can find neighborhoods of $f_0(x_0)$
and $F(x_0)$ respectively, over which the inverse functions are uniquely defined. We
denote by $f_0^{-1}$ and $F^{-1}$ respectively the inverses of $f_0$ and $\distfunc$ over
such a neighborhood of $f_0(x_0)$ and $F(x_0)$ respectively. The function $f_0$ is
assumed to be monotone on its whole support.

\begin{itemize}
  \item[\FincreasingText]\label{assumpt: Fincreasing}
        $\distfunc$ is locally continuous and strictly increasing in the sense that for all $x'$ in $C$,
        for all $\epsilon>0$, there exists $\gamma>0$ such that $\vert \distfunc^{-1}(u)-x'\vert >\gamma$
        for all $u$ such that $\vert u-\distfunc(x')\vert \geq \epsilon$.
  \item[\fdecreasingText]\label{assumpt: fdecreasing}
        $f_0$ is non-increasing, and $f_0$ is locally strictly decreasing in the sense that for all $x'$ in $C$,
        for all $\epsilon>0$, there exists $\gamma>0$ such that $\vert f_0(x')-f_0(y)\vert >\gamma$
        for all $y$ such that $\vert y-x'\vert \geq \epsilon$.
  \item[\fcontinuousText]\label{assumpt: fcontinuous}
        $f_0$ continuous in $C$.
\end{itemize}

Assumptions {\AssumpX}, {\AssumpC} and {\Fincreasing} ensure that $X_t$ visits infinitely
many times any small enough neighborhood of $x_0$ with probability 1, and guarantee that
$x_0$ is not at the boundary of the recurrent states. Assumptions {\AssumpX} and
{\AssumpC} and Lemma 3.2 in \cite{Tjostheim-2001} imply that \(T_n(C)\to\infty\) almost
surely.

\begin{theorem}\label{theo: cvfn2}
  Suppose that assumptions {\AssumpX}-{\fcontinuous} are satisfied.
  Then, as $n\to\infty$, one has
  \begin{equation}\label{eq: cvfndirect}
    \hat f_n(x_0)=f_0(x_0)+o_P(1),
  \end{equation}
  and
  \begin{equation}\label{eq: cvfninv}
    \hat f_n^{-1}(f_0(x_0))=x_0+o_P(1).
  \end{equation}
\end{theorem}
\section{Rates of convergence}\label{sec: rates}

To compute rates of convergence, we need stronger assumptions than for consistency. We
replace assumption {\AssumpX} for the following stronger version
\begin{itemize}
  \item[\AssumpXrateText]\label{assumpt: AssumpXrate}
        $\{X_t\}$ is a positive or \(\beta\)-null recurrent, aperiodic
        and irreducible Markov Chain whose kernel \(P(x,A)\) satisfies the minorization condition (\ref{eq: minorization}).
\end{itemize}

We replace, {\Fincreasing}, {\fdecreasing} and {\fcontinuous}, for the following slightly
more restrictive assumption

\begin{itemize}
  \item[\AssumpDeriveeText]\label{assumpt: AssumpDerivee}
        The function $f_0$ is non-increasing, the functions $f_0$ and
        $F_C$ are differentiable in \(C\), and the derivatives  $F_C'$ and
        $f_0'$ are bounded, in absolute value, above and away from zero
        in \(C\).
\end{itemize}

Let \(\lambda\) be the initial measure of \(\chain\). Our next hypothesis imposes some
control on the behavior of the chain in the first regenerative block.

\begin{itemize}
  \item[\AssumpInitialMeasureText]\label{assumpt: AssumpInitialMeasure}
        There exists a constant \(K\) and a neighborhood \(V\) of 0, such that
        \[\E_{\lambda}\left(\sum\limits_{t = 0}^{{\tau_{\check{\atom}} }} {\left( {{\I_C}\left\{ {{X_t} \leqslant {x_0} + \gamma } \right\} - {\I_C}\left\{ {{X_t} \leqslant {x_0} - \gamma } \right\}} \right)}\right) \leqslant K\gamma \quad \forall \gamma  \in V.\]
\end{itemize}

Assumption \AssumpInitialMeasure is satisfied if we assume that the initial measure of
the chain is the small measure used for the construction of the split chain (see equation
4.16c in \cite{Nummelin1984}). In the positive recurrent case, taking \(\lambda\) equal
to the unique invariant probability measure of the chain also satisfies
\AssumpInitialMeasure.

And finally, we need to control the number of times the chain visits \(C\) in a
regeneration block.

\begin{itemize}
  \item[\AssumpSizeBlockText]\label{assumpt: AssumpSizeBlock} \(\ell_C(\mathcal{B}_1)=\sum_{t\in \block_1}\I_C{\left\{ X_t \right\}}\) has finite second moment.
\end{itemize}

\begin{theorem}\label{theo: ratedirect}
  Assume that {\AssumpW}, {\AssumpC}, {\AssumpInteriorPoint},
  {\AssumpXrate}, {\AssumpDerivee}, {\AssumpInitialMeasure} and {\AssumpSizeBlock} hold.
  Then, as $n\to\infty$, one has
  \begin{equation}\label{eq: ratefn}
    \hat f_n(x_0)=f_0(x_0)+O_P(\asymptnumbreg^{-1/3}),
  \end{equation}
  with \(\asymptnumbreg\) as defined in \eqref{eq: definition_un}.
\end{theorem}

The rate $\asymptnumbreg$ comes from Lemmas \ref{lem: cvFnInvRate} and \ref{lem: Tn}, and
as it can be seen from Theorem \ref{th: t_n_general_results}, it is a deterministic
approximation of $\numreg$. Note that in the positive recurrent case, $\asymptnumbreg =
n$, hence we obtain the same rate $n^{-1/3}$ as in the i.i.d. case \cite[Chapter
2]{groeneboom2014nonparametric}. In the $\beta$-null recurrent case, however, the rate of
convergence is $n^{\beta/3}L^{1/3}\left(n\right)$ which is slower than the usual rate.
This is due to the null-recurrence of the chain because it takes longer for the process
to return to a neighborhood of the point $x_0$ and it is these points in the neighborhood
of $x_0$ which are used in nonparametric estimation.

\section{Localized Markov chains}\label{sec: localized_markov_chains}

Given the localized nature of our approach, in this section, we present some results that
are particularly useful in this scenario. These results are well known for positive
recurrent chains but are new in the null recurrent case. The detailed proofs of these
results can be found in Section \ref{sec: prooflocalized_markov_chains} of the
Supplementary material \cite{BertailDurotFernandezSupplements2023}.

The first result can be viewed as an extension of the Glivenko-Cantelli theorem to the
localized scenario.

\begin{lemma}
  \label{lem: cvFn}
  Assume that {\AssumpX} and {\AssumpC} hold. {Then, there exists a stationary $\sigma$-finite measure $\pi$, and $F$ defined by \eqref{eq: F},
      such that,
      \begin{equation}\label{eq: cvFnloc}
        \mathop {\sup }\limits_{y \in \R} \left| {{F_n}\left( y \right) - \distfunc\left( y \right)} \right| \to 0 \quad a.s.
      \end{equation}
      as $n\to\infty$.}
  If {\Fincreasing} is also satisfied, then, for all sufficiently
  small $\epsilon>0$, as $n\to\infty$ we have
  \begin{equation}\label{eq: cvFninvloc}
    \sup\limits_{\left|p - F\left(x_0\right)\right| \leqslant \epsilon } \left| {F_n^{ - 1}\left(p\right) - \distfunc^{-1}\left(p\right)} \right| \to 0\quad a.s.
  \end{equation}
\end{lemma}

Our next result (Lemma \ref{lem: covering+VC}), which is an extension of Lemma 2 in
\cite{bertail-portier:2019} to the localized \(\beta\)-null recurrent case, deals with
the properties of classes of functions defined over the regeneration blocks. Before
presenting the result, we need some machinery.

Recall that $E\subseteq\R$ denotes the state space of $X$. Define $\hat{E}= \cup_{k=1}
^{\infty}E^{k}$ (i.e. the set of finite subsets of \(E\)) and let the \textit{localized
occupation measure} $M_C$ be given by
\begin{align*}
  M_{C}(B, dy ) = \sum_{x\in B \cap C} \delta_{x}(y), \qquad\text{for every } B\in
  \hat{E}.
\end{align*}

The function that gives the size of the localized blocks is $\ell_C:
\hat{E}\to\mathbb{N}$
\begin{align*}
  \ell_C (B)= \int\, M_C(B,\mathrm{d} y),\qquad\text{ for every } B\in \hat{E}.
\end{align*}

Let $\hat{\SgE}$ denote the smallest $\sigma$-algebra formed by the elements of the
$\sigma$-algebras $\mathcal{E}^{k}$, $k\geq1$, where $\mathcal{E}^{k}$ stands for the
classical product $\sigma$-algebra. Let $\hat{Q}$ denote a probability measure on
$(\hat{E}, \hat{\mathcal{E}})$. If $B(\omega)$ is a random variable with distribution
$Q^{\prime }$, then $M_C (B(\omega),\mathrm{d} y)$ is a random measure, i.e., $M_C
(B(\omega), \mathrm{d} y )$ is a (counting) measure on $(E,\mathcal{E})$, almost surely,
and for every $A\in\mathcal{E}$, $M_C (B(\omega),A) = \int_{A} \, M_C(B(\omega
),\mathrm{d} y) $ is a measurable random variable (valued in $\mathbb{N}$). Henceforth
$\ell(B(\omega))\times\int f(y) \, M_C(B(\omega),\mathrm{d} y)$ is a random variable and,
provided that $\hat{Q}(\ell^{2} ) <\infty$, the map $Q_C$, defined by
\begin{align}
  \label{def:proba_measure}Q_C(A) = {E}_{\hat{Q}} \left(  \ell_C (B)
  \times\int_{A} \, M_C (B,\mathrm{d} y)\right)  / E_{\hat{Q}}( \ell_C^2)
  ,\qquad\text{for every } A\in\mathcal{E},
\end{align}

is a probability measure on $(E,\mathcal{E})$. The notation $E_{Q_C}$ stands for the
expectation with respect to the underlying measure $Q_C$. Introduce the following
notations: for any function $g: E\to \mathbb R$, let $\hat{g}_C:\hat{E} \to \mathbb R$ be
given by
\begin{equation}\label{eq:hat_function}
  \hat{g}_C \left( B \right) = \int {g\left( y \right){M_C}\left( {B,dy} \right)}  = \sum\limits_{x \in B \cap C} {g\left( x \right)}  = \sum\limits_{x \in B} {{g_C}\left( x \right)},
\end{equation}
and for any class $\mathcal{G}$ of real-valued functions defined on $E$,
denote the localized version of the sums on the blocks
by $\hat{G}_C  = \{ \hat{g}_C\,:\, g\in \mathcal G\} $.

Notice that, for any function \(g\),
\begin{equation}\label{eq:integral_q}
  {\E_{{Q_C}}}\left( g \right) = \frac{{{\E_{{\hat{Q}}}}\left( {{\ell_C}\left( B \right) \times \int {g\left( y \right){M_C}\left( {B,dy} \right)} } \right)}}{{{\E_{{\hat{Q} }}}\left( {\ell_C^2} \right)}} = \frac{{{\E_{{\hat{Q} }}}\left( {{\ell_C}\left( B \right)\hat{g}_C \left( B \right)} \right)}}{{{\E_{{\hat{Q} }}}\left( {\ell_C^2} \right)}}.
\end{equation}

\begin{lemma}
  \label{lem: covering+VC} Let $\hat{Q}$ be a probability measure
  on $(\hat{E}, \hat{\mathcal{E}})$ such that
  $0< \|\ell_C\|_{L^{2}(\hat{Q})}<\infty$ and $\mathcal{G}$ be a
  class of measurable real-valued functions defined on
  $(E,\mathcal E)$.
  Then we have, for every $0<\epsilon<\infty$,
  \[\mathcal{N}\left( {\epsilon {{\left\| {{\ell_C}} \right\|}_{{L_2}\left( {{\hat{Q} }} \right)}},\hat{\mathcal{G}}_C ,{L^2}\left( {{\hat{Q} }} \right)} \right) \leqslant \mathcal{N}\left( {\epsilon ,\mathcal{G},{L^2}\left( Q \right)} \right),\]
  where $Q$ is given in (\ref{def:proba_measure}). Moreover, if $\mathcal{G}$ belongs to
  the Vapnik–Chervonenkis (VC) class of functions with constant envelope $U$ and
  characteristic $(\textbf{C},v)$, then $\hat{\mathcal{G}}$ is VC with envelope $U \ell_C
  $ and characteristic $(\textbf{C},v)$.
\end{lemma}
\begin{remark}
  For a probability measure \(\mu\), and a class of functions \(\mathcal{H}\), the covering
  number \(\mathcal{N}\left( {\epsilon ,\mathcal{H},{L^r}\left( \mu \right)} \right)\) is the
  minimum number of \(L^r\left(\mu\right)\) \(\epsilon\)-balls needed to cover \(\mathcal{H}\).
  For more details about this concept and the VC class of functions, see \cite{Kosorok2008}.
\end{remark}

To put into perspective Lemma \ref{lem: covering+VC}, consider a class of bounded
functions \(\mathcal{G}\) that is VC with finite envelope. Lemma \ref{lem: covering+VC}
tells us that the class of unbounded functions \(\hat{\mathcal{G}}_C\) is also VC. If we
also have that {\AssumpSizeBlock} holds, then Theorem 2.5 in \cite{Kosorok2008} tells us
that \(\hat{\mathcal{G}}_C\) is a Donsker class. A reasoning like this is used in the
proof of the following result, which is a stronger version of Lemma \ref{lem: cvFn} under
assumptions {\AssumpXrate} and {\AssumpDerivee} and has some interest on its own.

\begin{lemma}
  \label{lem: cvFnInvRate}
  Assume that {\AssumpXrate}, {\AssumpDerivee}, {\AssumpC},
  {\AssumpInteriorPoint} and {\AssumpSizeBlock}
  hold. Then, for all sufficiently small  $\epsilon>0$ we have,
  \begin{equation}\label{eq: approxP}
    {T_n}(C)\mathop {\sup }\limits_{|y - x_0| \leq \epsilon } {\left| {F_n(y) - \distfunc(y)} \right|^2} = {O_p}\left( 1 \right)
  \end{equation}
  when \(n\) goes to \(+\infty\). If {\AssumpDerivee} is also satisfied, as \(n\to\infty\)
  we have
  \begin{equation}\label{eq: approxinvP}
    {T_n}(C)\mathop {\sup }\limits_{|p - F(x_0)| \leq \epsilon } {\left| {F_n^{ - 1}(p) - {\distfunc^{ - 1}}(p)} \right|^2} = {O_p}\left( 1 \right).
  \end{equation}
\end{lemma}
\section{Proofs} \label{sec: proofs}

In this section, we provide the proof of Theorems \ref{theo: cvfn2} and \ref{theo:
ratedirect}. These proofs make use of several intermediate lemmas, whose proofs can be
found in Sections \ref{sec: proofsbasic} and \ref{sec: proofsprelimrate}.

\subsection{Proof of Theorem \ref{theo: cvfn2}}\label{sec: basic}
Recall that we consider the piecewise constant and left-continuous
LSE $\hat f_n$, that is constant on every interval $(Y_{k-1},Y_{k}]$,
$k=2,\dots,m$ and also on $(-\infty,Y_{1}]$ and on $[Y_{m},\infty)$.
With $\delta>0$ fixed, we denote by $\Tn$  the number of times
the Markov Chain $\chain$ visits the set
$C:=[x_0-\delta,x_0+\delta]$ until time $n$:
\begin{equation}\label{eq: Tn}
  \Tn=\sum_{t=0}^n\I{\{X_t\in C\}}.
\end{equation}

Let $l_k = \sum_{t=1}^n \I_C{\left\{X_t\leq Y_k\right\}}$ for all $k\in\{1,\dots,m\}$ and
\(l_0=0\).

Our aim is to provide a characterization of $\hat f_n(x_0)$. Recall from \eqref{eq: Fn}
that the localized empirical distribution function \(F_n\) is defined as
\begin{equation*}
  F_n(y)=\frac 1{\Tn} \sum_{i=0}^{\Tn}{\I{\{ X_{\sigma_C\left(i\right)}\leq y \}}} = \frac 1{\Tn}\sum_{t=0}^n\I_C{\{X_t\leq y\}}
\end{equation*}
for \(y\in\R\). \(F_n\) is 0 on \(\left( { - \infty ,{Y_1}} \right)\),
so, with an arbitrary random variable \(Y_0<Y_1\) we have \(F_n(y)=F_n(Y_0)=0\)
for all \(y<Y_1\). Let  $\cal K$ be the set
\begin{equation}\label{eq: knotsties}
  \mathcal K:=\left\{F_n(Y_k),\ k=0,\dots,m\right\}
\end{equation}
and let $\Lambda_n$ be the continuous piecewise-linear
process on $[F_n(Y_0),F_n(Y_m)]$ with knots at the points in $\cal K$
and values
\begin{equation}\label{eq: Lambdaties}
  \Lambda_n\left(F_n(Y_k)\right)=\frac 1{\Tn}\sum_{t=0}^nZ_t\I_C{\left\{X_t\leq Y_k\right\}}
\end{equation}
at the knots. The characterization of $\hat f_n$ in Lemma \ref{lem: basic} involves
the least concave majorant of $\Lambda_n$. Note that we use $\Tn$ as a normalization
in the definitions of the processes $F_n$ and  $\Lambda_n$ since this choice ensures
that $F_n$ and $\Lambda_n$  converge to fixed functions, see
Lemma \ref{lem: cvFn}.

\begin{lemma}\label{lem: linear}
  For all \(y \in \left[ {F_n\left(Y_0\right),{F_n}\left( {{Y_m}} \right)} \right]\),
  \[{\Lambda _n}\left( y \right) = {L_n}\left( y \right) + {M_n}\left( y \right),\]
  where,
  \begin{equation}\label{eq: Ln1Extended}
    {L_n}\left( y \right) = \int\limits_{0}^y {f \circ F_{_n}^{ - 1}\left( u \right)du},
  \end{equation}

  and $M_n$ is a piece-wise linear processes with knots at $F_n(Y_k)$ for
  $k\in\{0,\dots,m\}$ such that
  \begin{equation*}
    M_n(F_n(Y_k))=\frac1{\Tn}\sum_{t=0}^nW_{t}\I_C{\{X_t\leq Y_k\}}.
  \end{equation*}

  Moreover, \(M_n\) can be written as

  \begin{equation}\label{eq: Mn1Extended}
    M_n(y)=\left\{
    \begin{array}{crl}
      0                & , & \textnormal{if}\; y=0  \\
      R_n^j(y) + M_n^j & , & \textnormal{otherwise}
    \end{array}
    \right.
  \end{equation}

  where,
  \begin{align}
    M_n^j               & = M_n(F_n(Y_j)) = \frac1{\Tn}\sum_{t=0}^n W_{t}\I_C{\{X_t\leq Y_j\}},\label{eq: Mn1}                                                          \\
    R_n^j\left(y\right) & = \frac{\sum\limits_{t = 0}^n {{W_t} \I_C{\{X_t = Y_{j+1}\}}}}{{{l_{j + 1}} - {l_j}}}\left( {y - F_n\left(Y_j\right)} \right)\label{eq: Rn1},
  \end{align}
  and \(j\) is such that \(Y_{j+1}=F_n^{-1}\left(y\right)\).
\end{lemma}

In the next lemma, we give an alternative characterization of the monotone nonparametric
LSE $\hat f_n$ at the observation points $Y_1,\dots,Y_m$.

\begin{lemma}\label{lem: basic}
  Let $C=[x-\delta,x+\delta]$ for some fixed $\delta>0$. Let $\hat\lambda_n$ be
  the left-hand slope of the least concave majorant of $\Lambda_n$. Then,
  \begin{equation}\label{eq: brunkties}
    \hat f_n(Y_k)=\hat\lambda_n \circ F_n(Y_k), \quad \forall k\in \{1,\dots,m\}.
  \end{equation}
  with probability 1 for \(n\) big enough.
\end{lemma}

We consider below the generalized inverse function of $\hat f_n$ since it has a more
tractable characterization than $\hat f_n$ itself. To this end, let us define precisely
the generalized inverses of all processes of interest. Since $\hat\lambda_n$ is a
non-increasing left-continuous step function on $(F_n(Y_0), F_n(Y_m)]$ that can have
jumps only at the points $F_n(Y_k)$, $k\in\{1,\dots,m\}$, we define its generalized
inverse $\hat U_n(a)$, for $a\in\R$, as the greatest $y\in(F_n(Y_0), F_n(Y_m)]$ that
satisfies $\hat\lambda_n(y)\geq a$, with the convention that the supremum of an empty set
is $F_n(Y_0)$. Then for every $a\in\R$ and $y\in(F_n(Y_0), F_n(Y_m)]$, one has
\begin{equation}\label{eq: inv}
  \hat \lambda_n(y)\geq a\mbox{  if and only if }\hat U_n(a)\geq y.
\end{equation}
Likewise, since $\hat f_n$ is a left-continuous non-increasing step function on $\R$ that can have jumps only at the observation times $Y_1<\dots<Y_m$, we define the generalized inverse $\hat f_n^{-1}(a)$, for $a\in\R$, as the greatest
$y\in[Y_0,Y_m]$ that satisfies $\hat f_n(y)\geq a$, with the convention that
the supremum of an empty set is $Y_0$. We then have
\begin{equation}\label{eq: invfn}
  \hat f_n(y)\geq a\mbox{  if and only if }\hat f_n^{-1}(a)\geq y
\end{equation}
for all $a\in\R$ and $y\in(Y_0,Y_m]$.  On the other hand, since $F_n$ is a right-continuous non-decreasing step function on $\R$ with range $[F_n(Y_0),F_n(Y_m)]$, we define the generalized inverse $F_n^{-1}(a)$, for $a\leq F_n(Y_m)$, as  the smallest $y\in[Y_0,Y_m]$ which satisfies $F_n(y)\geq a$. Note that the infimum is achieved for all $a\leq F_n(Y_m)$. We then have
\begin{equation}\label{eq: invFn}
  F_n(y)\geq a\mbox{  if and only if }F_n^{-1}(a)\leq y
\end{equation}
for all $a\leq F_n(Y_m)$ and $y\in[Y_0,Y_m]$, and thanks to Lemma \ref {lem: basic}  we have
\begin{equation}\label{eq: invm}
  \hat f_n^{-1}=F_n^{-1}\circ \hat U_n
\end{equation}
on $\R$.
Moreover, one can check that
\begin{equation}\label{eq: Un}
  \hat U_n(a)=\argmax_{p\in[F_n(Y_0),F_n(Y_m)]}\{\Lambda_n(p)-ap\},\mbox{ for all } a\in\R,
\end{equation}

where argmax denotes the greatest location of maximum (which is achieved on the set $\cal
K$ in \eqref{eq: knotsties}). Thus, the inverse process $\hat U_n$ is a location process
that is more tractable than $\hat f_n$ and $\hat \lambda_n$ themselves. A key idea in the
following proofs is to derive properties of $\hat U_n$ from its argmax characterization
(\ref{eq: Un}), then, to translate these properties to $\hat f_n^{-1}$ thanks to
(\ref{eq: invm}), and finally to translate them to $\hat f_n$ thanks to (\ref{eq:
invfn}).

To go from $\hat U_n$ to $\hat f_n^{-1}$ using \eqref{eq: invm} requires to approximate
$F_n^{-1}$ by a fixed function. Hence, in the sequel, we are concerned by the convergence
of the process $F_n$ given in \eqref{eq: Fn}, where $\delta>0$ is chosen sufficiently
small, and by the convergence of the corresponding inverse function $F_n^{-1}$.

It is stated in Lemma \ref{lem: cvFn} that under {\AssumpX} and {\AssumpC}, $F_n$
converges to a fixed distribution function $\distfunc$ that depends on $C$, hence on
$\delta$. If, moreover, $\distfunc$ is strictly increasing in $C$, then we can find a
neighborhood of $\distfunc(x_0)$ over which the (usual) inverse function $\distfunc^{-1}$
is uniquely defined, and $F_n^{-1}$ converges to $\distfunc^{-1}$.

In the following lemma, we show that $\distfunc(x_0)$ belongs to the domain of
$\Lambda_n$ with probability tending to one as $n\to\infty$. {
    \begin{lemma}
      \label{lem: Findomain}
      Assume that {\AssumpX}, {\AssumpC}, {\AssumpInteriorPoint} and {\Fincreasing} hold. Then, we can find
      $\epsilon>0$ such that the probability that
      $Y_1+\epsilon\leq x_0\leq Y_m-\epsilon$ tends to one as
      $n\to\infty$.
      Moreover, the probability that $F_n(Y_1)\leq \distfunc(x_0)\leq F_n(Y_m)$ tends
      to one as $n\to\infty$.
    \end{lemma}
  }

We will also need to control the noise $\{W_t\}$. The following lemma shows that the
noise is negligible under our assumptions.

\begin{lemma}\label{lem: Wnegligeable}
  Assume that {\AssumpX} and {\AssumpW} hold. Let \(\mathcal{F}_n=\sigma\left(\left\{ X_1, \ldots, X_n \right\}\right)\). Then,
  \[\sum\limits_{t = 0}^n {{W_t}{\I_C}{\left\{ {{X_t} = {A_n}} \right\}}}  = {o_P}\left( {\Tn} \right),\]
  and
  \[\mathop {\sup }\limits_{u > {A_n}} \left| {\sum\limits_{t = 0}^n {{W_t}{\I_C}\left\{ {{X_t} \in \left( {{A_n},u} \right]} \right\}} } \right| = {o_P}\left( {\Tn} \right).\]
  for any sequence of random variables \(A_n\), independent of the process $\{W_t\}$,
  that is adapted to the filtration \(\{ \mathcal{F}_n \}\).
\end{lemma}

With the above lemmas, we can prove convergence of $\hat U_n$ to $F\left(x_0\right)$
given by \eqref{eq: Un}.

\begin{lemma}\label{lem: cvUn}
  Suppose that assumptions {\AssumpX}-{\fcontinuous} are satisfied.
  Then, as $n\to\infty$, one has
  \begin{equation}\label{eq: cvUn}
    \hat U_n(f_0(x_0))=\distfunc(x_0)+o_P(1).
  \end{equation}
\end{lemma}

Now we proceed to the proof of (\ref{eq: cvfninv}). Fix $\epsilon>0$ arbitrarily small.
It follows from \eqref{eq: invm} and \eqref{eq: invFn} that
\begin{eqnarray*}
  \P\left(\hat f_n^{-1}(a)>x_0+\epsilon\right)&\leq&\P\left(F_n^{-1}\circ \hat U_n(a)>x_0+\epsilon\right)\\
  &\leq&\P\left( \hat U_n(a)\geq F_n(x_0+\epsilon)\right)\\
  &\leq&\P\left( \hat U_n(a)\geq \distfunc(x_0+\epsilon)-K_n\right),
\end{eqnarray*}
where
\begin{equation*}
  K_n=\sup_{\vert y-x_0\vert \leq\epsilon}\vert F_n(y)-\distfunc(y)\vert.
\end{equation*}
With $\nu:=\distfunc(x_0+\epsilon)-\distfunc(x_0)$, we obtain
\begin{eqnarray*}
  \P\left(\hat f_n^{-1}(a)>x_0+\epsilon\right)
  &\leq&\P\left( \hat U_n(a)\geq \distfunc(x_0)+\nu-K_n\right),
\end{eqnarray*}
and $\nu$ is strictly positive since $F$ is strictly increasing  in the neighborhood of $x_0$.
Hence, it follows from \eqref{eq: cvFnloc} that for sufficiently small $\epsilon>0$ one has
\begin{eqnarray*}
  \P\left(\hat f_n^{-1}(a)>x_0+\epsilon\right)
  &\leq&\P\left( \hat U_n(a)\geq \distfunc(x_0)+\nu/2\right)+o(1),
\end{eqnarray*}
so it follows from \eqref{eq: cvUn} that the probability
that $\hat f_n^{-1}(a)>x_0+\epsilon$ tends to zero as
$n\to\infty$. Similarly, the probability that $\hat f_n^{-1}(a)<x_0-\epsilon$
tends to zero as $n\to\infty$ so we conclude that the probability that
$\vert\hat f_n^{-1}(a)-x_0\vert>\epsilon$ tends to zero as $n\to\infty$.
This completes the proof of (\ref{eq: cvfninv}).

To prove (\ref{eq: cvfndirect}), fix $\epsilon>0$ sufficiently small so that $F$ and
$f_0$ are continuous and strictly increasing in the neighborhood of
$x':=f_0^{-1}(f_0(x_0)+\epsilon)$. Equation (\ref{eq: cvfninv}) shows that
\begin{equation}\label{eq: cvfneps}
  \hat f_n^{-1}(f_0(x_0)+\epsilon)=f_0^{-1}(f_0(x_0)+\epsilon)+o_P(1),
\end{equation}
as $n\to\infty$. Now, it follows from the switch relation \eqref{eq: inv} that
\begin{align}\notag\label{eq: invtodirect}
  \P\left(\hat f_n(x_0)>f_0(x_0)+\epsilon\right) & \leq \P\left(\hat f_n^{-1}(f_0(x_0)+\epsilon)\geq x\right)                               \\
                                                 & \leq\P\left(\hat f_n^{-1}(f_0(x_0)+\epsilon)\geq f_0^{-1}(f_0(x_0)+\epsilon)+\nu\right),
\end{align}
where $\nu:=x-f_0^{-1}(f_0(x_0)+\epsilon)>0$. It follows
from \eqref{eq: cvfneps} that the probability on the right-hand side tends to
zero as $n\to\infty$. Hence, the probability on the left-hand side tends
to zero as well as $n\to\infty$.

Similarly, the probability that $\hat f_n(x_0)<f_0(x_0)-\epsilon$ tends to zero as
$n\to\infty$ so we conclude that the probability that $\vert\hat
f_n(x_0)-f_0(x_0)\vert>\epsilon$ tends to zero as $n\to\infty$. This completes the proof
of Theorem \ref{theo: cvfn2}.

\subsection{Proof of Theorem \ref{theo: ratedirect}}\label{sec: prelimrate}

The proof of Theorem \ref{theo: ratedirect}, uses similar ideas as the ones used in the
proof of Theorem \ref{theo: cvfn2} but under stronger assumptions (and therefore using
stronger lemmas).

The first intermediate result is the following stronger version of Lemma \ref{lem:
Wnegligeable}.

\begin{lemma}\label{lem: Wnegligeablerate}
  Assume that {\AssumpW}, {\AssumpC}, {\AssumpInteriorPoint},
  {\AssumpXrate}, {\AssumpDerivee} and {\AssumpInitialMeasure} hold. Then, there exists $K>0$,
  $\gamma_0>0$ that do not depend on $n$ and \(N_{\gamma_0}\in\N\), such that for all
  $\gamma\in[0,\gamma_0]$ and \(n\geq N_{\gamma_0}\) one has
  \begin{align}
    \E_{\lambda}\left(\sup_{|y-x_0|\leq \gamma}\left\vert\sum_{t=0}^nW_{t}\left(\I_C{\{X_t\leq y\}}-\I_C{\{X_t\leq x_0\}}\right)\right\vert^2\right) & \leq K\asymptnumbreg\gamma  \label{eq: controlW}     \\
    \E_{\lambda}\left(\sup_{|y-x_0|\leq \gamma}\left\vert\sum_{t=0}^nW_{t}\I_C{\{X_t= y\}}\right\vert^2\right)                                       & \leq K\asymptnumbreg\gamma\label{eq: controlWfixedY}
  \end{align}
\end{lemma}

Then, we need to quantify how well we can approximate \(\Tn\) by \(\asymptnumbreg\).
\begin{lemma}\label{lem: Tn}
  Assume that {\AssumpXrate} and {\AssumpC} hold. Then we have
  \begin{itemize}
    \item[a)] As $n\to\infty$ we have
          \begin{equation*}
            \frac{\asymptnumbreg}{\Tn}=O_P(1).
          \end{equation*}

    \item[b)] Let \(\alpha\) and \(\eta\) be positive constants, then there positive
          exists constants \(N_\eta\), \(\underline{c}_{\eta}\) and
          \(\overline{c}_{\eta}\), such that
          \[\P\left( {{{\left( {\frac{{T_n\left( C \right)}}{{a\left( n \right)}}} \right)}^\alpha } \in \left[ {\underline{c}_{\eta},\overline{c}_{\eta}} \right]} \right) \geqslant 1 - \eta ,\quad \forall n \geqslant {N_\eta }.\]
  \end{itemize}

\end{lemma}

With the above lemmas (including Lemma \ref{lem: cvFnInvRate} and the ones used in
Section \ref{sec: basic}), we can obtain the rate of convergence of $\hat U_n$ given by
\eqref{eq: Un}.

\begin{lemma}\label{lem: rateinv}
  Assume that {\AssumpW}, {\AssumpC}, {\AssumpInteriorPoint},
  {\AssumpXrate}, {\AssumpDerivee}, {\AssumpInitialMeasure} and {\AssumpSizeBlock} hold. Then, as $n\to\infty$, one has
  \begin{equation}\label{eq: rateUn}
    \hat U_n(f_0(x_0))=\distfunc(x_0)+O_P\left(\asymptnumbreg^{-1/3}\right),
  \end{equation}
  and
  \begin{equation}\label{eq: rateinv}
    \hat f_n^{-1}(f_0(x_0))=x+O_P\left(\asymptnumbreg^{-1/3}\right).
  \end{equation}
\end{lemma}

Inspecting the proof of Lemma \ref{lem: rateinv}, one can see that the convergences in
\eqref{eq: rateUn} and \eqref {eq: rateinv} hold in a uniform sense in the neighborhood
of $x_0$. More precisely, there exists $\gamma>0$, independent on $n$, such that for all
$\eta>0$ we can find $K_1>0$ such that
\begin{equation*}
  \sup_{\vert a-f_0(x_0)\vert\leq\gamma}\P\left(\left\vert \hat U_n(a)-\distfunc\circ f_0^{-1}(a)\right\vert > K_1{\asymptnumbreg^{-1/3}}\right)\leq\eta
\end{equation*}
and
\begin{equation*}
  \sup_{\vert a-f_0(x_0)\vert\leq\gamma}\P\left(\left\vert \hat f_n^{-1}(a)-f_0^{-1}(a)\right\vert > K_1{\asymptnumbreg^{-1/3}}\right)\leq\eta.
\end{equation*}
Let $\epsilon= K_1{\asymptnumbreg^{-1/3}}$ where $K_1>0$ does not depend on $n$, and
recall \eqref{eq: invtodirect}
where $\nu=x-f_0^{-1}(f_0(x_0)+\epsilon)>0$. It follows from the assumption {\AssumpDerivee} that $f_0^{-1}$ has a derivative that is bounded in sup-norm away from zero in a neighborhood of $f_0(x_0)$. Hence, it follows from the Taylor expansion that there exists $K_2>0$ that depends only on $f_0$ such that $\nu\geq K_2\epsilon$, provided that $n$ is sufficiently large to ensure that $f_0(x_0)+\epsilon$ belongs to this neighborhood of $f_0(x_0)$. Hence,
\begin{eqnarray*}
  \P\left(\hat f_n(x_0)>f_0(x_0)+\epsilon\right)
  &\leq&\P\left(\hat f_n^{-1}(f_0(x_0)+\epsilon)\geq f_0^{-1}(f_0(x_0)+\epsilon)+K_2\epsilon\right).
  \\
  &\leq&\sup_{\vert a-f_0(x_0)\vert\leq\gamma}\P\left(\left\vert \hat f_n^{-1}(a)-f_0^{-1}(a)\right\vert > K_2K_1{\asymptnumbreg^{-1/3}}\right),
\end{eqnarray*}
provided that $n$ is sufficiently large to ensure that $f_0(x_0)+\epsilon$ belongs to the above neighborhood of $f_0(x_0)$, and that $\gamma\geq C{\asymptnumbreg^{-1/3}}$. For fixed $\eta>0$ one can choose $K_2>0$ such that the probability on the right-hand side of the previous display  is smaller than or equal to $\eta$ and therefore,
\begin{equation*}
  \lim_{n\to\infty}\P\left(\hat f_n(x_0)>f_0(x_0)+K_2{\asymptnumbreg^{-1/3}}\right)
  \leq\eta.
\end{equation*}
Similarly, for all fixed $\eta>0$, one can find $K_3$ that does not depend on $n$ such that
\begin{equation*}
  \lim_{n\to\infty}\P\left(\hat f_n(x_0)<f_0(x_0)-K_3{\asymptnumbreg^{-1/3}}\right)
  \leq\eta.
\end{equation*}
Hence, for all fixed $\eta>0$, there exists $K>0$ that independent of $n$ such that
\begin{equation*}
  \lim_{n\to\infty}\P\left(\vert\hat f_n(x_0)-f_0(x_0)\vert>K{\asymptnumbreg^{-1/3}}\right)
  \leq\eta.
\end{equation*}This completes the proof of Theorem \ref{theo: ratedirect}.

\section{Technical proofs for Section \ref{sec: basic}}\label{sec: proofsbasic}

\par {\bf Proof of Lemma \ref{lem: linear}.}
Combining \eqref{eq: Lambdaties} and \eqref{eq: model} yields
\begin{equation*}
  \Lambda_n(F_n(Y_k))
  =\frac1{\Tn}\sum_{t=0}^nf_0(X_{t})\I_C{\{X_t\leq Y_k\}} +\frac1{\Tn}\sum_{t=0}^n W_{t}\I_C{\{X_t\leq Y_k\}}.
\end{equation*}

The first term on the right-hand side of the previous display can be rewritten as
follows:
\begin{align*}
  \frac{1}{\Tn}\sum_{t=0}^nf_0(X_{t})\I_C{\{X_t\leq Y_k\}} & =\frac1{\Tn}\sum_{j=1}^mf_0(Y_l)(l_j-l_{j-1})\I_C{\{Y_j\leq Y_k\}} \\
                                                           & =\sum_{j=1}^k\int _{l_{j-1}/\Tn}^{l_j/\Tn} f_0\circ F_n^{-1}(u)du,
\end{align*}

using that $F_n^{-1}(u)=Y_j$ for all $u\in(l_{j-1}/\Tn,l_j/\Tn]$. Hence, for all $k$ in
$\{0,\dots,m\}$
\begin{equation}\label{eq: Lambdanties}
  \Lambda_n(F_n(Y_k))= \int_{0}^{l_k/\Tn}f_0\circ F_n^{-1}(u) du +\frac1{\Tn}\sum_{t=0}^nW_{t}\I_C{\{X_t\leq Y_k\}}.
\end{equation}
Combining \eqref{eq: Lambdanties} with the piece-wise linearity of $\Lambda_n$ yields
\[{\Lambda _n}\left( {{F_n}\left( {{Y_k}} \right)} \right) = {L_n}\left( {{F_n}\left( {{Y_k}} \right)} \right) + {M_n}\left( {{F_n}\left( {{Y_k}} \right)} \right),\]
where $L_n$ and $M_n$ are piece-wise linear processes with knots at $F_n(Y_k)$ for $k$
  in $\{0,\dots,m\}$ and such that
\begin{equation*}
  L_n(F_n(Y_k))=\int_{0}^{l_k/\Tn}f_0\circ F_n^{-1}(u) du
\end{equation*}
and
\begin{equation*}
  M_n(F_n(Y_k))=\frac1{\Tn}\sum_{t=0}^nW_{t}\I_C{\{X_t\leq Y_k\}}.
\end{equation*}

In order to ease the notation, we will write $F_n^i=F_n(Y_i)$, $L_n^i = {L_n}\left(
{{F_n}\left( {{Y_i}} \right)} \right)$ and $M_n^i = {M_n}\left( {{F_n}\left( {{Y_i}}
\right)} \right)$. Let $y\in \left( {{F_n}\left( {{Y_0}} \right),{F_n}\left( {{Y_m}}
\right)} \right]$, take \(j\) such that \(Y_{j+1}=F_n^{-1}\left(y\right)\), then
${F_n}\left( {{Y_j}} \right) < y \leqslant {F_n}\left( {{Y_{j + 1}}} \right)$. With this
notation,
\begin{align*}
  {L_n}\left( y \right) & = \frac{{L_n^{j + 1} - L_n^j}}{{F_n^{j + 1} - F_n^j}}\left( {y - F_n^j} \right) + L_n^j, \\
  {M_n}\left( y \right) & = \frac{{M_n^{j + 1} - M_n^j}}{{F_n^{j + 1} - F_n^j}}\left( {y - F_n^j} \right) + M_n^j.
\end{align*}
Notice that
\begin{align*}
  L_n^{j + 1} - L_n^j & = \int\limits_{\frac{{{l_j}}}{{\Tn}}}^{\frac{{{l_{j + 1}}}}{{\Tn}}} {f_0 \circ {F_n^{-1}}\left( u \right)du}  = \frac{{{l_{j + 1}} - {l_j}}}{{\Tn}}f\left( {{Y_{j + 1}}} \right), \\
  F_n^{j + 1} - F_n^j & = \frac{{{l_{j + 1}} - {l_j}}}{{\Tn}},
\end{align*}
therefore,
\begin{equation*}
  {L_n}\left( y \right) = f_0\left( {{Y_{j + 1}}} \right)\left( {y - F_n^j} \right) + L_n^j = \int\limits_{\frac{{{l_j}}}{{\Tn}}}^y {f_0 \circ F_{_n}^{ - 1}\left( u \right)du}  + L_n^j = \int\limits_{0}^y {f_0 \circ F_{_n}^{ - 1}\left( u \right)du},
\end{equation*}
which proves (\ref{eq: Ln1Extended}).

For \(M_n\) we have,
\begin{equation*}
  M_n^{j + 1} - M_n^j = \frac{1}{{\Tn}}\sum\limits_{t = 0}^n {{W_t} \I_C{\{X_t = Y_{j+1}\}}},
\end{equation*}
then,
\begin{equation*}
  {M_n}\left( y \right) = \frac{\sum\limits_{t = 1}^n {{W_t} \I_C{\{X_t = Y_{j+1}\}}}}{{{l_{j + 1}} - {l_j}}}\left( {y - F_n^j} \right) + M_n^j = R_n^j(y)+M_n^j.
\end{equation*}
and this completes the proof.
\par {\bf Proof of Lemma \ref{lem: basic}.}
By definition, with $l_0=0$,
and $l_k = \sum_{t=0}^n \I_C{\left\{X_t\leq Y_k\right\}}$ for all $k\in\{1,\dots,m\}$,
we have $F_n(Y_k)=al_k$ for all $k\in\{0,\dots,m\}$, where $a=1/\Tn$
and does not depend on $k$.
Moreover,
\begin{equation*}
  \Lambda_n\left(F_n(Y_k)\right)= a\sum_{t=0}^nZ_t\I_C{\left\{X_t\leq Y_k\right\}}
\end{equation*}
Since $\hat f_n(Y_k)$ is the left-hand slope at $l_k$ of the least concave majorant  of the set of points in \eqref{eq: CSD}, the equality in \eqref{eq: brunkties} follows from Lemma 2.1 in \cite{durot2003distance}.
\par{\bf Proof of Lemma \ref{lem: Findomain}.}
The first assertion follows from
Assumption {\AssumpInteriorPoint} and the second immediately follows from the
first one by \eqref{eq: cvFnloc} combined with the strict monotonicity
of $F$ in \(C\).
\par{\bf Proof of Lemma \ref{lem: Wnegligeable}.}
Let \(\mathcal{F}_n=\sigma\left(\left\{ X_0, \ldots, X_n \right\}\right)\) be sigma
algebra generated by the chain \( \left\{ X_t \right\} \) up to time \(n\). Denote by
\(\P_{\mathcal{F}_n}\) the probability conditioned to \(\mathcal{F}_n\). Take \(\epsilon>0\).

By Chebyshev's inequality, we have
\[{\P_{{\mathcal{F}_n}}}\left( {\left| {\frac{{\sum\limits_{t=0}^n {{W_t}{\I_C}\left\{ {{X_t} = {A_n}} \right\}} }}{{{T_n}\left( C \right)}}} \right| > \epsilon} \right) \leqslant \frac{{{\sigma ^2}\sum\limits_{t = 0}^n {{\I_C}\left\{ {{X_t} = {A_n}} \right\}} }}{{{\epsilon^2}{T_n}^2\left( C \right)}} \leqslant \frac{{{\sigma ^2}}}{{{\epsilon^2}{T_n}\left( C \right)}},\]
which implies the first part of the Lemma because \(\Tn\to\infty\) with probability 1.

For the second part, let \(\gamma_n(u)\) be the number of times the chain visits
\({\left( {{A_n},u} \right]}\cap C\) up to time \(n\) and \({A_n}\left( u \right) =
\left\{ {t \leqslant n:{X_t} \in \left( {{A_n},u} \right]}\cap C \right\}=\left\{
a_1,\ldots, a_{\gamma_n(u)} \right\}\) the times of those visits. Using that
\(\gamma_n=\mathop {\sup }\limits_{u > {A_n}} {\gamma _n}\left( u \right) \leqslant
{T_n}\left( C \right)\) and Kolmogorov's inequality (Th 3.1.6, pp 122 in \cite{Gut2013})
we obtain,
\begin{align*}
  {\P_{{\mathcal{F}_n}}}\left( {\mathop {\sup }\limits_{u > {A_n}} \left| {\frac{{\sum\limits_{t = 0}^n {{W_t}{\I_C}\left\{ {{X_t} \in \left( {{A_n},u} \right]} \right\}} }}{{{T_n}\left( C \right)}}} \right| > \epsilon } \right) & = {\P_{{\mathcal{F}_n}}}\left( {\mathop {\sup }\limits_{u > {A_n}} \left| {\sum\limits_{i = 1}^{{\gamma _n}\left( u \right)} {\frac{{{W_{{t_{{a_i}}}}}}}{{{T_n}\left( C \right)}}} } \right| > \epsilon } \right)  \\
                                                                                                                                                                                                                                     & \leqslant {\P_{{\mathcal{F}_n}}}\left( {\mathop {\sup }\limits_{1 \leqslant k \leqslant \gamma_n} \left| {\sum\limits_{i = 1}^k {\frac{{{{W}_{t_{a_i}}}}}{{{T_n}\left( C \right)}}} } \right| > \epsilon } \right) \\
                                                                                                                                                                                                                                     & \leqslant \frac{{{\sigma ^2}}}{{{\epsilon^2}{T_n}\left( C \right)}}.
\end{align*}

which by the same argument as before, implies the second part of the Lemma.
\par{\bf Proof of Lemma \ref{lem: cvUn}.}
In the sequel, we set $a=f_0(x_0)$. We begin with the proof of \eqref{eq: cvUn}.

Fix $\epsilon>0$ arbitrarily, and let $\nu>0$ and $\gamma>0$ be such that $\vert
\distfunc^{-1}(u)-x_0\vert >\nu$ for all $u$ such that $\vert u-\distfunc(x_0)\vert \geq
\epsilon/2$, and $\vert f_0(x_0)-f_0(y)\vert >\gamma$ for all $y$ such that $\vert
y-x_0\vert \geq \nu/2$. Note that existence of $\nu$ and $\gamma$ is ensured by
assumptions {\Fincreasing} and {\fdecreasing}.

By Lemma \ref{lem: Findomain}, we can assume without loss of generality that
$\distfunc(x_0)$ belongs to the domain $[F_n(Y_1),F_n(Y_m)]$ of $\Lambda_n$, since this
occurs with probability tending to one. Therefore, we can find $j(x_0)$ such that
${Y_{j(x_0)}} = {F_n}^{ - 1}\left( {\distfunc\left( {{x_0}} \right)} \right)$. It follows
from the characterization in \eqref{eq: Un} that the event $E_n^1:=\{\hat U_n(a)>
\distfunc(x_0)+\epsilon\}$ is contained in the event that there exists $p\in\cal K$ such
that $p>\distfunc(x_0)+\epsilon$ and
\begin{equation*}
  {\Lambda _n}\left( p \right) - ap \geqslant {\Lambda _n}\left( {\distfunc\left( {{x_0}} \right)} \right) - a\distfunc\left( {{x_0}} \right),
\end{equation*}
where we recall that $a=f_0\left(x_0\right)$.

By Lemma \ref{lem: linear}, $E_n^1$ is contained in the event that there exists $p\in\cal
K$ such that $p>\distfunc(x_0)+\epsilon$ and
\begin{equation}\label{eq: maximalCondition}
  L_n(p)+M_n(p)-ap\geq L_n(\distfunc(x_0))+M_n(\distfunc(x_0))-a\distfunc(x_0)
\end{equation}
Using \eqref{eq: Ln1Extended} in \eqref{eq: maximalCondition} we obtain that $E_n^1$ is
contained in the event that there exists $p\in\cal K$ such that
$p>\distfunc(x_0)+\epsilon$ and
\begin{equation*}
  \int_{t_0/\Tn}^{p}f_0\circ F_n^{-1}(u) du+S_n-ap\geq \int_{t_0/\Tn}^{\distfunc(x_0)}f_0\circ F_n^{-1}(u) du-a\distfunc(x_0),
\end{equation*}
where
\begin{eqnarray*}
  {S_n} = \mathop {\sup }\limits_{p > \distfunc(x_0) + \epsilon ,\;p \in \mathcal{K} } \left\{ {{M_n}(p) - {M_n}\left( {\distfunc\left( {{x_0}} \right)} \right)} \right\}.
\end{eqnarray*}
Let \(j\) and \(k\) such that \(Y_{j+1}=F_n^{-1}\left(F\left(x_0\right)\right)\) and $p=F_n(Y_k)$. By equation \eqref{eq: Mn1Extended} we have ${M_n}\left( p \right) - {M_n}\left( {\distfunc\left( {{x_0}} \right)} \right) = M_n^k - M_n^j - R_n^j(\distfunc(x_0))$, therefore,
\begin{align*}
  {S_n} & = \mathop {\sup }_{\begin{subarray}{l}
                               p > \distfunc\left( {{x_0}} \right) + \epsilon \\
                               p \in \mathcal{K}
                             \end{subarray}}  {\left\{ {M_n^k - M_n^j} \right\} - R_n^j\left( {\distfunc\left( {{x_0}} \right)} \right)}\nonumber \\
  \begin{split}
     & \leq \mathop {\sup }_{\begin{subarray}{l}
                               p > \distfunc\left( {{x_0}} \right) + \epsilon \\
                               p \in \mathcal{K}
                             \end{subarray}}{\left\vert\frac{1}{\Tn}\sum_{t=0}^nW_{t}\left(\I_C{\{X_t\leq F_n^{-1}(p)\}}-\I_C{\{X_t\leq F_n^{-1}(\distfunc(x_0))\}}\right)\right\vert} \\
     & \qquad + \left| {R_n^j\left( {F_n\left(Y_{j + 1}\right)} \right)} \right|\nonumber                                                                                              \\
  \end{split}
  \\[2ex]
  \begin{split}
     & \leq \mathop {\sup }_{\begin{subarray}{l}
                               p > \distfunc\left( {{x_0}} \right) + \epsilon \\
                               p \in \mathcal{K}
                             \end{subarray}}{\left\vert\frac{1}{\Tn}\sum_{t=0}^nW_{t}\I_C{\left\{X_t\in \left( {F_n^{ - 1}\left( {\distfunc\left( {{x_0}} \right)} \right);F_n^{ - 1}(p)} \right] \right\}}\right\vert} \\
     & \qquad + \frac{\left\vert\sum\limits_{t = 0}^n {{W_t} \I_C{\left\{X_t = F_n^{ - 1}\left( {\distfunc\left( {{x_0}} \right)} \right)\right\}}}\right\vert}{\Tn}.
  \end{split}
\end{align*}
Hence,
\begin{equation*}
  \begin{split}
    \Tn S_n & \leq \mathop {\sup }_{\begin{subarray}{l}
                                      p > \distfunc\left( {{x_0}} \right) + \epsilon \\
                                      p \in \mathcal{K}
                                    \end{subarray}}{\left\vert \sum_{t=0}^nW_{t}\I_C{\left\{X_t\in \left( {F_n^{ - 1}\left( {\distfunc\left( {{x_0}} \right)} \right);F_n^{ - 1}(p)} \right] \right\}} \right\vert} \\
            & \qquad+ \left\vert\sum\limits_{t = 0}^n {{W_t} \I_C{\left\{X_t = F_n^{ - 1}\left( {\distfunc\left( {{x_0}} \right)} \right)\right\}}}\right\vert.
  \end{split}
\end{equation*}
Therefore, the event $E_n^1$ is contained in the event that there exists
$p>\distfunc(x_0)+\epsilon$ such that
\begin{equation*}
  \int_{\distfunc(x_0)}^{p}f_0\circ F_n^{-1}(u) du+S_n\geq a(p-\distfunc(x_0)).
\end{equation*}
Now, let $E_n^2$ be the event that
\begin{equation*}
  \sup_{\vert u-\distfunc(x_0)\vert\leq\epsilon}\vert F_n^{-1}(u)-F^{-1}(u)\vert\leq \eta
\end{equation*}
where $\eta\in(0,\nu/4)$ is such that $\vert f_0(y)-f_0(x_0)\vert\leq \gamma/2 $ for all $y$ such that $\vert x_0-y\vert \leq \eta$. Note that the existence of $\eta$ is ensured by assumption \fcontinuous.
Then, it follows from the monotonicity of $f_0$ and $F_n$ that on $E_n^2$,
\begin{eqnarray*}
  \int_{\distfunc(x_0)}^{p}f_0\circ F_n^{-1}(u) du&\leq&
  \int_{\distfunc(x_0)}^{\distfunc(x_0)+\epsilon/2}f_0( F^{-1}(u)-\eta) du+\int_{\distfunc(x_0)+\epsilon/2}^{p}f_0( F_n^{-1}(\distfunc(x_0)+\epsilon/2)) du.
\end{eqnarray*}
Hence, it follows from the definitions of $\eta$, $\nu$ and $\gamma$ that on $E_n^2$,
\begin{eqnarray*}
  \int_{\distfunc(x_0)}^{p}f_0\circ F_n^{-1}(u) du&\leq&
  \frac {\epsilon}{2}f_0(x_0)+\frac{\gamma\epsilon}4 +(p-\distfunc(x_0)-\epsilon/2)f_0(F^{-1}(\distfunc(x_0)+\epsilon/2)-\eta)\\
  &\leq&
  \frac {\epsilon}{2}f_0(x_0)+\frac{\gamma\epsilon}4 +(p-\distfunc(x_0)-\epsilon/2)f_0(x_0+\nu/2)\\
  &\leq&
  \frac {\epsilon}{2}f_0(x_0)+\frac{\gamma\epsilon}4 +(p-\distfunc(x_0)-\epsilon/2)(f_0(x_0)-\gamma).
\end{eqnarray*}
This implies that on $E_n^2$,
\begin{eqnarray*}
  \int_{\distfunc(x_0)}^{p}f_0\circ F_n^{-1}(u) du
  &\leq&
  a(p-\distfunc(x_0))-(p-\distfunc(x_0)-3\epsilon/4)\gamma\\
  &\leq&
  a(p-\distfunc(x_0))-\epsilon\gamma/4
\end{eqnarray*}
for all $p>\distfunc(x_0)+\epsilon$.
Hence, the event $E_n^1\cap E_n^2$ is contained in the event $\{S_n\geq \epsilon\gamma/4\}$. Now, on $E_n^2$, for all $p>\distfunc(x_0)+\epsilon$  we have
\begin{eqnarray*}
  F_n^{-1}(p)&\geq& F_n^{-1}(\distfunc(x_0)+\epsilon)\\
  &\geq& F^{-1}(\distfunc(x_0)+\epsilon)-\eta\\
  &\geq&x+\nu-\eta\\
  &\geq&F_n^{-1}(\distfunc(x_0))+\nu-2\eta\\
  &\geq&F_n^{-1}(\distfunc(x_0))+\nu/2,
\end{eqnarray*}
since $\nu>4\eta$. Therefore,
\begin{align*}
  \Tn S_n\leq & \sup_{u>F_n^{-1}(\distfunc(x_0))+\nu/2}\left\vert\sum_{t=0}^nW_{t}\I_C{\{X_t\in(F_n^{-1}(\distfunc(x_0)),  u]\}}\right\vert + \\&+  \left\vert\sum\limits_{t = 0}^n {{W_t} \I_C{\{X_t = F_n^{ - 1}\left( {\distfunc\left( {{x_0}} \right)} \right)\}}}\right\vert.
\end{align*}
Hence, it follows from Lemma \ref{lem: Wnegligeable} that $S_n$ converges in probability to zero as $n\to\infty$, so that the probability of the event $\{S_n\geq \epsilon\gamma/4\}$ tends to zero as $n\to\infty$.
It follows from Lemma \ref{lem: cvFn} that for $\epsilon$ sufficiently small, the probability of the event $E_n^2$ tends to one as $n\to\infty$, so we conclude that the probability of $E_n^1$ tends to zero as $n\to\infty$. Similarly, the probability of the event $\{\hat U_n(a)< \distfunc(x_0)-\epsilon\}$ tends to zero as $n\to\infty$, so that
\begin{equation*}
  \lim_{n\to\infty}\P(\vert\hat U_n(a)- \distfunc(x_0)\vert>\epsilon)=0
\end{equation*}
for all $\epsilon>0$. This completes the proof of \eqref{eq: cvUn}.

\section{Technical proofs for Section \ref{sec: prelimrate}}\label{sec: proofsprelimrate}

\par{\bf Proof of Lemma \ref{lem: Wnegligeablerate}.}
Let \(\mathcal{F}_n=\sigma\left(\left\{ X_0, \ldots, X_n \right\}\right)\) be sigma
algebra generated by the chain \(\chain\) up to time \(n\). Denote by
\(\E_{\mathcal{F}_n}\) the expected value conditioned to \(\mathcal{F}_n\).
Take \(0<\gamma\leq\delta\) and define \({I_0 } = \left[ {{x_0} - \gamma ,{x_0}} \right]\),
\({I_1 } = \left[ {{x_0} ,{x_0}+ \gamma} \right]\) and
\[\begin{gathered}
    {S_0\left(\gamma\right) } = \mathop {\sup }\limits_{y \in {I_0 }} \left| {\sum\limits_{t=0}^n {{W_t}\left( {{\I}\left\{ {{X_t} \leqslant y} \right\} - {\I}\left\{ {{X_t} \leqslant {x_0}} \right\}} \right)} } \right|^2 \hfill \\
    {S_1\left(\gamma\right) } = \mathop {\sup }\limits_{y \in {I_1 }} \left| {\sum\limits_{t=0}^n {{W_t}\left( {{\I}\left\{ {{X_t} \leqslant y} \right\} - {\I}\left\{ {{X_t} \leqslant {x_0}} \right\}} \right)} } \right|^2 \hfill \\
  \end{gathered} \]
then,
\begin{align}
  S\left(\gamma\right) & =\mathop {\sup }\limits_{\left| {y - {x_0}} \right|\leq\gamma } \left| {\sum\limits_{t=0}^n {{W_t}\left( {{\I}\left\{ {{X_t} \leqslant y} \right\} - {\I}\left\{ {{X_t} \leqslant {x_0}} \right\}} \right)} } \right|^2=\max{\Bigl(S_0\left(\gamma\right), S_1\left(\gamma\right)\Bigr)},\nonumber \\
                       & \leq S_0\left(\gamma\right)+S_1\left(\gamma\right)\label{eq: supremum-decomposition}
\end{align}

Following the notation of section \ref{sec: markov_theory}, let
\[\alpha_n^{(0)}\left(\gamma\right) = \mathop {\sup }\limits_{y \in {I_0 }} {T_n}\left( {\left[ {y,{x_0}} \right]} \right)\quad ,\quad \alpha_n^{(1)}\left(\gamma\right)  = \mathop {\sup }\limits_{y \in {I_1 }} {T_n}\left( {\left[ {{x_0},y} \right]} \right),\]
with this notation, \({S_0} = \mathop {\sup }\limits_{y \in {I_0}} {\left|
  {\sum\limits_{i=1}^{{T_n}\left( {\left[ {y,{x_0}} \right]} \right)}
  {{W_{{\sigma_{[y,x_0)}(i)}}}} } \right|^2}\) and \(S_1 = \mathop {\sup }\limits_{y \in
  {I_1 }} \left| {\sum\limits_{{i=1}}^{{T_n}\left( {\left[ {{x_0},y} \right]} \right)}
  {{W_{{{\sigma_{[x_0,y]}(i)}}}}} } \right|^2\).

By Doob's maximal inequality (Th 10.9.4 in \cite{Gut2013}), we have, for \(j=0,1\),
\begin{align}\label{eq: bound_expected_value_s_j}
  {\E_{{\mathcal{F}_n}}}S_j\left(\gamma\right) & \leqslant 4{\E_{{\mathcal{F}_n}}}{\left( {\sum\limits_{i = 1}^{{\alpha^{(j)}_n}} {{W_{{t_i}}}} } \right)^2}=4\sigma^2\alpha^{(j)}_n\left(\gamma\right) \nonumber \\
                                               & \leq 4{\sigma ^2}{T_n}\left( {\left[ {{x_0} - \gamma ,{x_0} + \gamma } \right]} \right)  \nonumber                                                               \\
                                               & \leq 4{\sigma ^2}\sum\limits_{t=0}^n {\left( {{\I}\left\{ {{X_t} \leqslant {x_0} + \gamma } \right\} - {\I}\left\{ {{X_t} < {x_0} - \gamma } \right\}} \right)}.
\end{align}

Therefore, by (\ref{eq: supremum-decomposition}) and (\ref{eq: bound_expected_value_s_j})
\begin{equation}\label{eq: bound_expected_value}
  {\E_{{\mathcal{F}_n}}}S\left(\gamma\right) \leqslant 8{\sigma ^2}\sum\limits_{t=0}^n {\Bigl( {{\I}\left\{ {{X_t} \leqslant {x_0} + \gamma } \right\} - {\I}\left\{ {{X_t} < {x_0} - \gamma } \right\}} \Bigr)}.
\end{equation}

Define,
\begin{itemize}
  \item \(h\left( {y,\gamma } \right) = \I{\left\{ {y \in \left[ {{x_0} - \gamma ,{x_0} + \gamma } \right]} \right\}}\),
  \item \(h\left( {{\block_j},\gamma } \right) = \begin{cases}
          \sum\limits_{t=0}^{{\tau _\alpha }} {h\left( {{X_t},\gamma } \right)}                                  & ,\quad j = 0 \hfill         \\
          \sum\limits_{t = {\tau _A}(j) + 1}^{{\tau _A}\left( {j + 1} \right)} {h\left( {{X_t},\gamma } \right)} & ,\quad j \geqslant 1 \hfill \\
        \end{cases}\)
  \item \({Z_n}\left( \gamma  \right) = \sum\limits_{t=0}^n {h\left( {{X_t},\gamma } \right)} \)
  \item \(\ell\left( {{\block_j}} \right) = \begin{cases}
          {\tau_{\atom} }                           & ,\quad j = 0 \hfill         \\
          {\tau_{\atom}}(j + 1) - {\tau_{\atom}}(j) & ,\quad j \geqslant 1 \hfill \\
        \end{cases}\)
  \item \(\widetilde T\left( n \right) = \min \left\{ {k:\sum\limits_{i = 0}^k {\ell\left( {{\block_j}} \right)}  \geqslant n} \right\}\).
  \item \({\mathcal{G}_k} = \sigma \left( {\left\{ {\left( {h\left( {{\block_j},\gamma } \right),\ell\left( {{\block_j}} \right)} \right)} \right\}_{j = 0}^k} \right)\) for \(k\geq 0\).
\end{itemize}

By the Strong Markov property, \(\left\{ {\left( {h\left( {{\block_j},\gamma }
\right),\ell\left( {{\block_j}} \right)} \right)} \right\}_{j = 1}^{ + \infty }\) is an
i.i.d. sequence which is independent of \({\left( {h\left( {{\block_0},\gamma }
\right),\ell\left( {{\block_0}} \right)} \right)}\) (and, therefore, of the initial
measure \(\lambda\)). For \(n\) fixed, the random variable \(\widetilde T\left( n
\right)\) is a stopping time for the sequence \(\left\{ {\left( {h\left(
{{\block_j},\gamma } \right),\ell\left( {{\block_j}} \right)} \right)} \right\}_{j = 0}^{
+ \infty }\), in effect
\begin{align*}
  \left\{ {\widetilde T\left( n \right) = 0} \right\} & = \left\{ {\ell\left( {{\block_0}} \right) \geqslant n} \right\} \in {\mathcal{G}_0},                                                                                                                                                                        \\
  \left\{ {\widetilde T\left( n \right) = k} \right\} & = \bigcap\limits_{j = 0}^{k - 1} {\left\{ {\sum\limits_{i = 0}^j {\ell\left( {{\block_i}} \right)}  < n} \right\}}  \bigcap \left\{ {\sum\limits_{i = 0}^k {\ell\left( {{\block_j}} \right)}  \geqslant n} \right\} \in {\mathcal{G}_k}\quad\forall k\geq 1.
\end{align*}

For each \(n\) and \(\gamma\) we have that
\begin{align}\label{eq: z_n_decomposition}
  {Z_n}\left( \gamma  \right) & = \sum\limits_{t=0}^{{\tau _\alpha }} {h\left( {{X_t},\gamma } \right)}  + \sum\limits_{j = 1}^{T\left( n \right)} {h\left( {{\block_j},\gamma } \right)}  + \sum\limits_{t = {t_\alpha }\left( {T\left( n \right)} \right) + 1}^n {h\left( {{X_t},\gamma } \right)}\nonumber \\
                              & \leqslant {h\left( {{\block_0},\gamma } \right)} + \sum\limits_{j = 1}^{\widetilde T\left( n \right)} {h\left( {{\block_j},\gamma } \right)}.
\end{align}
where the last inequality is justified by the fact that, \(T\left( n \right) \leqslant \widetilde T\left( n \right)\)
and \(h\left(y,\gamma\right)\) is a nonnegative function. Because \(\ell\left( {{\block_j}} \right)\geq 1\) for all \(j\), we have that,
\[\sum\limits_{j = 1}^{\widetilde T\left( n \right)} {h\left( {{\block_j},\gamma } \right)}  = \sum\limits_{j = 1}^{n} {h\left( {{\block_j},\gamma } \right)\I\left\{ {\widetilde T\left( n \right)\geqslant j} \right\}},\]
then,
\begin{equation}\label{eq: z_n_expectation}
  \E\left( {\sum\limits_{j = 1}^{\widetilde T\left( n \right)} {h\left( {{\block_j},\gamma } \right)} } \right) = \sum\limits_{j = 1}^{n} {\E\left( {h\left( {{\block_j},\gamma } \right)\I\left\{ {\widetilde T\left( n \right)\geqslant j} \right\}} \right)}.
\end{equation}

For each \(j\) we have,
\[\E_{\lambda}\left( {h\left( {{\block_j},\gamma } \right)\I\left\{ {\widetilde T\left( n \right) \geqslant j} \right\}} \right) = \E_{\lambda}\left( {\E\left( {h\left( {{\block_j},\gamma } \right)\I\left\{ {\widetilde T\left( n \right) \geqslant j} \right\}\left| {{\mathcal{G}_{j - 1}}} \right.} \right)} \right)\]
Notice that \(\I\left\{ {\widetilde T\left( n \right) \geqslant j} \right\} = 1 -
  \I\left\{ {\widetilde T\left( n \right) \leqslant j - 1} \right\} \in {\mathcal{G}_{j -
  1}}\) and \({h\left( {{\block_j},\gamma } \right)}\) is independent of
  \({\mathcal{G}_{j - 1}}\), therefore,
\begin{align*}
  \E_{\lambda}\left( {h\left( {{\block_j},\gamma } \right)\I\left\{ {\widetilde T\left( n \right) \geqslant j} \right\}} \right) = \E_{\lambda}\left( {\I\left\{ {\widetilde T\left( n \right) \geqslant j} \right\}} \right)\E\left( {h\left( {{\block_j},\gamma } \right)} \right).
\end{align*}

Plugging this into equation \eqref{eq: z_n_expectation} we get,
\[\E_{\lambda}\left( {\sum\limits_{j = 1}^{\widetilde T\left( n \right)} {h\left( {{\block_j},\gamma } \right)} } \right) = \sum\limits_{j = 1}^{n} {\E\left( {h\left( {{\block_j},\gamma } \right)} \right)\P_{\lambda}\left( {\widetilde T\left( n \right) \geqslant j} \right)}  \leq \E\left( {h\left( {{\block_1},\gamma } \right)} \right)\E_{\lambda}\widetilde T\left( n \right).\]

Then, by taking expectation in \eqref{eq: z_n_decomposition} we obtain
\begin{align}\label{eq: z_n_expectation_bound}
  \E_{\lambda}{Z_n}\left( \gamma  \right) & \leqslant \E_{\lambda} h\left( {{\block_0},\gamma } \right) + \E\left( {h\left( {{\block_1},\gamma } \right)} \right)\E_{\lambda}\widetilde T\left( n \right)\nonumber \\
                                          & \leq \E_{\lambda} h\left( {{\block_0},\gamma } \right) + \E\left( {h\left( {{\block_1},\gamma } \right)} \right)\E_{\lambda}\left(T\left(n\right)+1\right).
\end{align}
By Theorem \ref{th:suplement:t_n_general_results} and the fact that \(\distfunc\) is Lipschitz we can find \(K_1\)
independent of \(\gamma\) such that,
\begin{align}
  \E\left( {h\left( {{\block_1},\gamma } \right)} \right) & = \int {h\left( {t,\gamma } \right)d\pi \left( t \right) = K_\pi \pi \left( C \right)\left( {{\distfunc}\left( {{x_0} + \gamma } \right) - {\distfunc}\left( {{x_0} - \gamma } \right)} \right)}\nonumber \\
                                                          & \leqslant K_1\gamma.\label{eq: expectation_lipschitz_bound}
\end{align}
If \(\chain\) is positive recurrent, by Theorem \ref{th:suplement:t_n_general_results},
\(\frac{\numreg}{\asymptnumbreg}\) converges almost surely to a positive constant \(K_2>0\).
Moreover,
\(\frac{T\left(n\right)}{\asymptnumbreg}\leq 1\) therefore, by the
Dominated Convergence Theorem we obtain that \(\E_{\lambda} T\left(n\right)\sim \frac{\asymptnumbreg}{K_2}\).
If \(\chain\) is \(\beta\)-null recurrent, by Lemma 3.3 in \cite{Tjostheim-2001},
\(\E_{\lambda} T\left( n \right) \sim \frac{{\asymptnumbreg}}{{\Gamma \left( {1 + \beta } \right)}}\),
hence, for both positive and \(\beta\)-null recurrent chains, we can find \(K_2\)
and \(N\), both independent of \(\gamma\), such that
\(\E_{\lambda} T\left( n \right) \leqslant {K_2}\asymptnumbreg\) for all \(n\geq N\).
Using this with \eqref{eq: z_n_expectation_bound} and \eqref{eq: expectation_lipschitz_bound} we get,
\begin{equation}\label{eq: z_n_final_expectation_bound}
  \frac{{\E_{\lambda}{Z_n}\left( \gamma  \right)}}{{\asymptnumbreg\gamma }} \leqslant \frac{{\E_{\lambda} h\left( {{\block_0},\gamma } \right)}}{{\asymptnumbreg\gamma }} + {K_1}{K_2}\quad \forall n \geqslant N,\;\forall \gamma  \in \left( {0,\delta } \right].
\end{equation}

Combining \eqref{eq: z_n_final_expectation_bound} with assumption {\AssumpInitialMeasure}
and the fact that \(Z_n\left(0\right)\equiv 0\) we obtain that there exist positive
constants \(K_3\) and \(\gamma_0\) such that
\begin{equation*}
  \E_{\lambda}{Z_n}\left( \gamma  \right) \leqslant \asymptnumbreg\gamma \quad \forall n \geqslant N,\;\forall \gamma  \in \left( {0,\gamma_0 } \right].
\end{equation*}
Equation \eqref{eq: controlW} now follows after taking expectation
in (\ref{eq: bound_expected_value}). The proof of \eqref{eq: controlWfixedY} follows the same reasoning, but using
\[{S_j}\left( \gamma  \right) = \mathop {\sup }\limits_{y \in {I_j}} {\left| {\sum\limits_{t=0}^n {{W_t}\left( {\I_C\left\{ {{X_t} = y} \right\}} \right)} } \right|^2}{\text{ }}.\]
\par{\bf Proof of Lemma \ref{lem: Tn}.}
a) If \(\chain\) is positive recurrent, Theorem \ref{th:suplement:t_n_general_results} implies that
there exists a positive constant \(K\) such that
\(\frac{T_n\left(C\right)}{\asymptnumbreg}\) converges almost surely to \(K\pi(C)\),
which is not zero by {\AssumpC}.

On the other hand, if \(\chain\) is \(\beta\)-null recurrent, Theorem
\ref{th:suplement:t_n_general_results} and Slutsky's Theorem implies that there exists a
constant \(K>0\) such that \(\frac{{{T_n}\left( C \right)}}{{\asymptnumbreg}}\) converges
in distribution to \(KM_{\beta}(1)\) where \(M_{\beta}(1)\) denotes a Mittag-Leffler
distribution with parameter \(\beta\). This distribution is continuous and strictly
positive with probability 1, then, by the Continuous Mapping Theorem,
\(\frac{{\asymptnumbreg}}{{{T_n}\left( C \right)}}\) converges in distribution to a
multiple of \(\frac{1}{{{M_\beta }}}\), therefore, \(\frac{{\asymptnumbreg}}{{{T_n}\left(
C \right)}}\) is bounded in probability by Theorem 2.4 in \cite{van2000asymptotic}.

b) Let \(\chain\) be positive recurrent, then, we can find \(N_\eta\) such that
\[\P\left( {\left| {{{\left( {\frac{{{T_n}\left( C \right)}}{{\asymptnumbreg}}} \right)}^\alpha } - K^{\alpha}\pi {{\left( C \right)}^{\alpha}}} \right| \leqslant {{\left( {\frac{{K\pi \left( C \right)}}{2}} \right)}^\alpha }} \right) \geqslant 1 - \eta,\quad \forall n \geqslant {N_\eta }.\]
hence,
\[\P\left( {{{\left( {\frac{{{T_n}\left( C \right)}}{{\asymptnumbreg}}} \right)}^\alpha } \in \left[ {\frac{{{K^{\alpha}\pi {{\left( C \right)}^\alpha }}}}{2},\frac{{3K^{\alpha}\pi {{\left( C \right)}^\alpha }}}{2}} \right]} \right)\geq 1 - \eta,\quad \forall n\geq{N_\eta }.\]

Now let \(\chain\) be \(\beta\)-null recurrent. Let \(Z={\left( {K{M_\beta }\left( 1
\right)} \right)^{\alpha}}\), This random variable is continuous and positive, therefore,
we can find positive constants \(\underline{c}_\eta\) and \(\overline{c}_\eta\) such that
\begin{equation}\label{eq: prob_mittag_leffler}
  \P\left( {Z \in \left[ {\underline{c}_\eta ,\overline{c}_\eta} \right]} \right) \geqslant 1 - \frac{\eta }{2}.
\end{equation}

By the Continuous Mapping Theorem, \({{{\left( {\frac{{{T_n}\left( C
\right)}}{{\asymptnumbreg}}} \right)}^{\alpha}}}\) converges in distribution to \(Z\),
therefore, we can find \(N_\eta\in\N\) such that
\begin{equation}\label{eq: convergence_mittag_leffler}
  \left| {\P\left( {{{\left( {\frac{{\Tn}}{{\asymptnumbreg}}} \right)}^\alpha } \in \left[ {{\underline{c}_\eta },{{\overline{c}}_\eta }} \right]} \right) - \P\left( {Z \in \left[ {{\underline{c}_\eta },{{\overline{c}}_\eta }} \right]} \right)} \right| \leq \frac{\eta }{2},\quad \forall n \geqslant {N_\eta },
\end{equation}

Combining \eqref{eq: prob_mittag_leffler} and \eqref{eq: convergence_mittag_leffler} we
obtain that
\begin{equation}\label{eq: approx_in_proba_t_n}
  \P\left( {{{\left( {\frac{{\Tn}}{{\asymptnumbreg}}} \right)}^\alpha } \in \left[ {{\underline{c}_\eta },{{\overline{c}}_\eta }} \right]} \right) \geqslant 1 - \eta ,\quad \forall n \geqslant {N_\eta }.
\end{equation}
\par{\bf Proof of Lemma \ref{lem: rateinv}.}
Fix $\epsilon\in(0,1)$ small enough so that $F'$ and $\vert f_0'\vert$ are bounded from above and away from zero on $[F^{-1}(\distfunc(x_0)-2\epsilon),F^{-1}(\distfunc(x_0)+2\epsilon)]$, see the assumption \AssumpDerivee.
Then, the proper inverse functions of $F$ and $f_0$ are well defined on $[\distfunc(x_0)-2\epsilon,\distfunc(x_0)+2\epsilon]$ and
\begin{equation*}
  [f_0\circ F^{-1}(\distfunc(x_0)-2\epsilon),f_0\circ F^{-1}(\distfunc(x_0)+2\epsilon)]
\end{equation*}
respectively. We denote the inverses on that intervals by $F^{-1}$ and $f_0^{-1}$ respectively.  Let
\begin{equation}\label{eq: Unloc}
  U_n(a)=\argmax_{\vert p-\distfunc(x_0)\vert\leq\epsilon}\{\Lambda_n(p)-ap\}
\end{equation}
where $a=f_0(x_0)$ and where the supremum is restricted to $p\in[F_n(Y_0),F_n(Y_m)]$.
We will show below that
\begin{equation}\label{eq: rateUnloc}
  U_n(a)=\distfunc(x_0)+O_P({\asymptnumbreg^{-1/3}}),
\end{equation}
as $n\to\infty$. Combining  \eqref{eq: Un} to Lemma \ref{lem: cvUn} ensures that $\hat U_n(a)$ coincides with $U_n(a)$ with a probability that tends to one as $n\to\infty$, so   \eqref{eq: rateUn} follows from \eqref{eq: rateUnloc}.

We turn to the proof of \eqref{eq: rateUnloc}. Fix $\eta>0$ arbitrarily and let
\begin{equation}\label{eq: gamma}
  \gamma_n=K_0 {\asymptnumbreg^{-1/3}} \end{equation}
for some $K_0\geq 1$ sufficiently large so that
\begin{equation}\label{eq: condgamma}
  \gamma_n\geq \frac{1}{\sqrt{\asymptnumbreg}} .\end{equation}
Then, by part ii) of Lemma \ref{lem: Tn}, we can find positive constants
$\underline c_\eta$, $\bar c_\eta$ and \(N_\eta\)
such that
\begin{equation}\label{eq: control_t_n_c}
  \P\left(T_n(C)^{2/3}\gamma_n \asymptnumbreg^{-1/3}\in[K_0\underline c_\eta,K_0\bar c_\eta]\right)\geq 1- \eta/2 \quad\forall n\geq N_\eta,
\end{equation}
Let \(\underline c = K_0 \underline c_\eta\) and \(\bar c = K_0 \overline c_\eta\).
It follows from \eqref{eq: approxinvP} that for sufficiently small
$\epsilon>0$, we can find $K_1>0$ such that
\begin{equation*}
  \P\left( T_n(C)\sup_{\vert p-\distfunc(x_0)\vert\leq 2\epsilon}\vert F_n^{-1}(p)-F^{-1}(p)\vert^2\leq K_1\right)\geq 1- \eta/2
\end{equation*}
for all $n$. Hence for \(n\geq\N_\eta\),
\begin{equation*}
  \P(\mathcal E_n)\geq 1-\eta,\end{equation*}
where $\mathcal E_n$ denotes the intersection of the events
\begin{equation}\label{eq : En1}
  T_n(C)^{2/3}\gamma_n \asymptnumbreg^{-1/3}\in[\underline c,\bar c]
\end{equation}
and
\begin{equation}\label{eq : En2}
  T_n(C)\sup_{\vert p-\distfunc(x_0)\vert\leq 2\epsilon}\vert F_n^{-1}(p)-F^{-1}(p)\vert^2\leq K_1.
\end{equation}
Combining equations (\ref{eq : En1}) and (\ref{eq : En2}), we obtain that, in \(\mathcal{E}_n\),
\begin{equation}\label{eq: bound_sup_e_n}
  \mathop {\sup }\limits_{\left| {p - {\distfunc}\left( {{x_0}} \right)} \right| \leqslant 2\varepsilon } {\left| {F_n^{ - 1}\left( p \right) - {F^{ - 1}}\left( p \right)} \right|^2} \leqslant {K_2}a{\left( n \right)^{ - 1}}
\end{equation}
where \(K_2={K_1}{\left( {\frac{{{K_0}}}{{\underline c }}} \right)^{3/2}}\) is
independent of \(n\) and \(K_0\).

By Lemma \ref{lem: Findomain}, we can assume without loss of generality that
$\distfunc(x_0)$ belongs to $[F_n(Y_0),F_n(Y_m)]$, since this occurs with probability
that tends to one. Hence, by (\ref{eq: Unloc}), the event $\{|U_n(a)-\distfunc(x_0)|\geq
\gamma_n\}$ is contained in the event that there exists $p\in[F_n(Y_0),F_n(Y_m)]$ with
$\vert p-\distfunc(x_0)\vert\leq\epsilon$, $|p-\distfunc(x_0)|\geq \gamma_n$ and
\begin{equation}
  \label{eq: argmax}
  \Lambda_n(p)-ap\geq \Lambda_n(\distfunc(x_0))-a\distfunc(x_0).
\end{equation}
Obviously, the probability is equal to zero if $\gamma_n>\epsilon$ so we assume in the sequel that $\gamma_n\leq\epsilon$. For all $p\in[\distfunc(x_0)-\epsilon,\distfunc(x_0)+\epsilon]$ define
\begin{equation*}
  \Lambda\left(p\right)=\int_{\distfunc\left({x_0}\right)}^p f_0\circ F^{-1}(u)du.
\end{equation*}

Let $c>0$ such that $\vert f_0'\vert/F'>2c$ on the interval
$[F^{-1}(\distfunc(x_0)-2\epsilon),F^{-1}(\distfunc(x_0)+2\epsilon)]$. Since
$\Lambda'(\distfunc(x_0))=a$ and $\Lambda''=f_0'\circ F^{-1}/F'\circ F^{-1}$, it then
follows from Taylor's expansion that
\begin{equation*}
  \Lambda(p)-\Lambda(\distfunc(x_0))\leq (p-\distfunc(x_0))a-c(p-\distfunc(x_0))^2
\end{equation*}
for all $p\in[\distfunc(x_0)-\epsilon,\distfunc(x_0)+\epsilon]$ and therefore, \eqref{eq: argmax} implies that
\begin{equation*}
  \Delta_n(p)-\Delta_n(\distfunc(x_0))
  -c(p-\distfunc(x_0))^2\geq0
\end{equation*}
for all such $p$'s, where we set $\Delta_n:=\Lambda_n-\Lambda$.
Hence, for all \(n\geq N_\eta\),
\begin{eqnarray}\label{eq1}\nonumber
  &&\P\left(\left|U_n\left(a\right)-\distfunc\left(x_0\right)\right|\geq \gamma_n\right)\\
  &&\qquad\leq\eta+ \P\left(\sup_{|p-\distfunc(x_0)|\in[\gamma_n,\epsilon]}\{\Delta_n(p)-\Delta_n(\distfunc(x_0))
  -c(p-\distfunc(x_0))^2\}\geq0\mbox{ and }{\mathcal E_n}\right)\nonumber\\
  &&\qquad\leq \eta+\sum_{j} \P\left(
  \sup_{|u|\in[\gamma_n2^{j},\gamma_n  2^{j+1}]}\{\Delta_n(\distfunc(x_0)+u)-\Delta_n(\distfunc(x_0))\}
  \geq c(\gamma_n 2^{j})^2  \mbox{ and }{\mathcal E_n}\right)\nonumber\\
  &&\qquad {\leq \eta+\sum_{j} \P\left(
    \sup_{|u|\leq\gamma_n  2^{j+1}}\vert \Delta_n(\distfunc(x_0)+u)-\Delta_n(\distfunc(x_0))\vert
    \geq c(\gamma_n 2^{j})^2  \mbox{ and }\mathcal E_n\right)}
\end{eqnarray}
where the sums are taken over all integers $j\geq0$ such that  $ \gamma_n2^{j}\leq\epsilon$.
Recall that we have \eqref{eq: Lambdanties}  for all $k\in\{0,\dots,m\}$.
Since  $\Lambda_n$ is
piecewise-linear with knots at $F_n(Y_0),\dots,F_n(Y_m)$, by (\ref{eq: Ln1Extended}) and (\ref{eq: Mn1Extended}) we get that
for every $j$ in the above sum,
\begin{eqnarray} \label{eq: decomposeMn}\notag
  &&\sup_{|u|\leq\gamma_n  2^{j+1}}|\Delta_n(\distfunc(x_0)+u)-\Delta_n(\distfunc(x_0))|\\ \notag
  &&\qquad  \leq
  \mathop {\sup }\limits_{|u| \leqslant {\gamma _n}{2^{j + 1}}} \left| {\int\limits_{\distfunc\left(x_0\right)}^{{\distfunc}\left( {{x_0}} \right) + u} {\left( {{f_0} \circ F_n^{ - 1}\left( y \right) - {f_0} \circ {F^{ - 1}}\left( y \right)} \right)dy} } \right|\\
  &&\qquad\qquad +\mathop {\sup }\limits_{|u| \leqslant {\gamma _n}{2^{j + 1}}} \left| {{M_n}\left( {{\distfunc}\left( {{x_0}} \right) + u} \right) - {M_n}\left( {{\distfunc}\left( {{x_0}} \right)} \right)} \right|.
\end{eqnarray}
Moreover,   $\vert f_0'\vert $ is bounded above on $[F^{-1}(\distfunc(x_0)-2\epsilon),F^{-1}(\distfunc(x_0)+2\epsilon)]$, so we obtain that
for every $j$ with $ \gamma_n2^{j}\leq\epsilon$, the first term on the right-hand side of \eqref{eq: decomposeMn} satisfies
\begin{eqnarray*}
  && \sup_{|u|\leq \gamma_n 2^{j+1}}\left\vert\int_{\distfunc(x_0)}^{\distfunc(x_0)+u}\left(f_0\circ F_n^{-1}(p)-f_0\circ F^{-1}(p)\right) dp\right\vert\\
  &&\qquad \leq\int_{\distfunc(x_0)-\gamma_n 2^{j+1}}^{\distfunc(x_0)+\gamma_n 2^{j+1}}\left\vert f_0\circ F_n^{-1}(p)-f_0\circ F^{-1}(p) \right\vert dp\\
  &&\qquad \leq
  K_3\gamma_n 2^{j}\sup_{\vert p-\distfunc(x_0)\vert\leq 2\epsilon}\left|F_n^{-1}(p)-F^{-1}(p)\right|,
\end{eqnarray*}
for some $K_3>0$ that does not depend on $n$. Hence,
it follows from the previous display and (\ref{eq: bound_sup_e_n}) that
\begin{eqnarray*}
  && \E\left(\sup_{|u|\leq \gamma_n 2^{j+1}}\left\vert\int_{\distfunc(x_0)}^{\distfunc(x_0)+u}(f_0\circ F_n^{-1}(p)-f_0\circ F^{-1}(p) dp\right\vert^2\I(\mathcal E_n)\right)\\
  &&\qquad \leq
  K_3^2\gamma_n^2 2^{2j}\E\left( \sup_{\vert p-\distfunc(x_0)\vert\leq 2\epsilon}|F_n^{-1}(p)-F^{-1}(p)|^2\I({\mathcal E_n})\right)\\
  &&\qquad \leq
  K_3^2\gamma _n^2{2^{2j}}{K_2}{\asymptnumbreg}^{ - 1}.
\end{eqnarray*}
Hence,  taking $K_4=K_3^2K_2$ we get that for all $j$ with $\gamma_n2^j\leq\epsilon\leq 1$.\begin{equation}\label{eq: controlesp}
  \E\left(\sup_{|u|\leq \gamma_n 2^{j+1}}\left\vert\int_{\distfunc(x_0)}^{\distfunc(x_0)+u}(f_0\circ F_n^{-1}(p)-f_0\circ F^{-1}(p) dp\right\vert^2\I(\mathcal E_n)\right)\leq
  {K_4}{\gamma_n}{2^j}\asymptnumbreg^{-1}.
\end{equation}
By equations (\ref{eq: Mn1Extended}) and (\ref{eq: Mn1}) in Lemma \ref{lem: linear},
the second term on the right-hand side of \eqref{eq: decomposeMn} satisfies,
\begin{equation}\label{eq: second_term}
  \mathop {\sup }\limits_{|u|\leq{\gamma _n}{2^{j + 1}}} \left| {{M_n}\left( {{\distfunc}\left( {{x_0}} \right) + u} \right) - {M_n}\left( {{\distfunc}\left( {{x_0}} \right)} \right)} \right| \leqslant {I^{n,j}_1} + {I^{n,j}_2},
\end{equation}
where \(I_1^{n,j}\) and \(I_2^{n,j}\) are given by
\begin{align*}
  I^{n,j}_1   & =\frac{1}{T_n(C)}\sup_{|u|\leq \gamma_n 2^{j+1}}\left\vert\sum_{t=0}^nW_{t}\Bigl(\I_C{\{X_t\leq F_n^{-1}(\distfunc(x_0)+u)\}}-\I_C{\{X_t\leq F_n^{-1}(\distfunc(x_0))\}}\Bigl)\right\vert,                                               \\
  {I^{n,j}_2} & = \frac{2}{{{T_n}  \left( C \right)}}\sup_{|u|\leq \gamma_n 2^{j+1}}\left\vert\sum\limits_{t = 0}^n {{W_t}\Bigl( {\I_C{\left\{ {{X_t} = F_n^{ - 1}\left( {{\distfunc}\left( {{x_0}} \right) + u} \right)} \right\}}} \Bigl)}\right\vert.
\end{align*}
For \(I^{n,j}_1\), it follows from the triangle inequality that
\begin{equation*}
  I^{n,j}_1 \leq  \frac{2}{T_n(C)}
  \sup_{|u|\leq \gamma_n 2^{j+1}}\left\vert\sum_{t=0}^nW_{t}\left(\I_C{\{X_t\leq F_n^{-1}(\distfunc(x_0)+u)\}}-\I_C{\{X_t\leq x_0\}}\right)\right\vert.
\end{equation*}
Combining (\ref{eq: bound_sup_e_n}) and the fact that
\(F^{-1}\) is Lipschitz in $[\distfunc(x_0)-2\epsilon, \distfunc(x_0)+2\epsilon]$ we can
find \(K_5=\max{\left(\sqrt{K_2}, \sup\left(F^{-1}\right) \right)}\) independent of \(n\)  such that, on $\mathcal E_n$,
\begin{equation*}\sup_{\vert p-\distfunc(x_0)\vert\leq 2\epsilon}\vert F_n^{-1}(p)-F^{-1}(p)\vert\leq \frac{K_5}{\sqrt{\asymptnumbreg}}
\end{equation*}
and $\vert F^{-1}(\distfunc(x_0)+u)-x\vert =\vert F^{-1}(\distfunc(x_0)+u)-F^{-1}(\distfunc(x_0))\vert\leq K_5\vert u\vert/2$ for all $u$ with $\vert u\vert\leq 2\epsilon$. Hence, on \(\mathcal{E}_n\)
\begin{align*}
  I^{n,j}_1 & \leq \frac{2}{T_n(C)}\sup_{|y-x_0|\leq K_5\gamma_n 2^j+K_5/\sqrt{\asymptnumbreg}}\left\vert\sum_{t=0}^nW_{t}\left(\I_C{\{X_t\leq y\}}-\I_C{\{X_t\leq x_0\}}\right)\right\vert,                             \\
  I^{n,j}_2 & \leq  \frac{2}{{{T_n}  \left( C \right)}}\sup_{|y-x_0|\leq K_5\gamma_n 2^j+K_5/\sqrt{\asymptnumbreg}}\left\vert\sum\limits_{t = 0}^n {{W_t}\Bigl( {\I_C{\left\{ {{X_t} = y} \right\}}} \Bigl)}\right\vert.
\end{align*}

It follows from \eqref{eq: condgamma} that $\gamma_n 2^j\geq \gamma_n\geq
1/\sqrt{\asymptnumbreg}$ for all $j\geq 0$, then, on $\mathcal E_n$
\begin{align*}
  I^{n,j}_1 & \leq \frac{2}{T_n(C)}\sup_{|y-x_0|\leq 2K_5\gamma_n 2^j}\left\vert\sum_{t=0}^nW_{t}\left(\I_C{\{X_t\leq y\}}-\I_C{\{X_t\leq x_0\}}\right)\right\vert,                             \\
  I^{n,j}_2 & \leq  \frac{2}{{{T_n}  \left( C \right)}}\sup_{|y-x_0|\leq 2K_5\gamma_n 2^j}\left\vert\sum\limits_{t = 0}^n {{W_t}\Bigl( {\I_C{\left\{ {{X_t} = y} \right\}}} \Bigl)}\right\vert.
\end{align*}

By Lemma \ref{lem: Wnegligeablerate}, we conclude that there exists $K_6>0$ and
\(N_\eta^\prime\) such that, for \(n\geq N_\eta^\prime\)
\begin{equation}\label{eq: bound_I_1_I_2}
  \E\Bigl(\left(I^{n,j}_1+I^{n,j}_2\right)^2\I(\mathcal E_n)\Bigr)\leq K_6\gamma_n 2^{j}\asymptnumbreg^{-1}
\end{equation}
Combining \eqref{eq: decomposeMn}, \eqref{eq: controlesp}, \eqref{eq: second_term} and \eqref{eq: bound_I_1_I_2}, we conclude
that there exists $K_7>0$, independent of \(n\) and \(K_0\), such that for
all \(n\geq N_\eta^\prime\) and $j\geq 0$ where $ \gamma_n2^{j}\leq\epsilon$, one has
\begin{eqnarray*}
  \E\left(\sup_{|u|\leq \gamma_n 2^{j+1}}|\Delta_n(\distfunc(x_0)+u)-\Delta_n(\distfunc(x_0))|^2{ \I(\mathcal E_n)}\right)\leq
  {K_7\gamma_n 2^{j}\asymptnumbreg^{-1}}.
\end{eqnarray*}
Combining this with \eqref{eq1} and the Markov inequality, we conclude that
there exist $K_8>0$ and \(N_\eta^{\prime\prime}\), that do not depend on $n$
nor \(K_0\), such that, for all \(n\geq N_\eta^{\prime\prime}\),
\begin{eqnarray*}
  \P\left(|U_n(a)-\distfunc(x_0)|\geq \gamma_n\right)&\leq&
  \eta+K_8\sum_{k\geq 0} \frac{\gamma_n 2^{j}\asymptnumbreg^{-1}}{(\gamma_n 2^{j})^4}\\
  &\leq&
  \eta+K_8\gamma_n^{-3}\asymptnumbreg^{-1}\sum_{j\geq 0} 2^{-3j}.
\end{eqnarray*}
The sum on the last line is finite, so there exists $K>0$, independent of \(n\)
and \(K_0\), such that for \(n\) bigger than \(N_\eta^{\prime\prime}\)
\begin{equation}\label{eq: final_prob_bound}
  \P\left(|U_n(a)-\distfunc(x_0)|\geq \gamma_n\right)\leq\eta+K\gamma_n^{-3}\asymptnumbreg^{-1}=\eta+\frac{K}{K_0^3}.
\end{equation}
The above probability can be made smaller than $2\eta$ by setting
\eqref{eq: gamma} for some sufficiently large $K_0$ independent of
$n$. This proves \eqref{eq: rateUnloc} and  completes
the proof of \eqref{eq: rateUn}.\\
Now, we turn to the proof of \eqref{eq: rateinv}.  It follows from \eqref{eq: invm} combined to
\eqref{eq: cvUn}
and Lemma \ref{lem: cvFnInvRate} that
\begin{equation*}
  \hat f_n^{-1}(f_0(x_0))=F^{-1}\circ \hat U_n(f_0(x_0))+T_n(C)^{-1/2}O_P\left(1\right).
\end{equation*}
Hence, by Lemma \ref{lem: Tn} we have
\begin{equation*}
  \hat f_n^{-1}(f_0(x_0))=F^{-1}\circ \hat U_n(f_0(x_0))+O_P\left(\asymptnumbreg^{-1/2}\right).
\end{equation*}
Now, it follows from the assumption {\AssumpDerivee} that $F^{-1}$ has a bounded derivative in the neighborhood of $\distfunc(x_0)$, to which $\hat U_n(f_0(x_0))$ belongs with probability that tends to one. Hence, it follows from Taylor's expansion that
\begin{align*}
  \hat f_n^{-1}(f_0(x_0)) & =F^{-1}\circ \distfunc(x_0)+O\left(\vert \hat U_n(f_0(x_0))-\distfunc(x_0)\vert\right)+O_P\left(\asymptnumbreg^{-1/2}\right) \\
                          & =x+O_P({\asymptnumbreg^{-1/3}})+O_P\left(\asymptnumbreg^{-1/2}\right),
\end{align*}
where we used \eqref{eq: rateUn} for the last equality. This proves \eqref{eq: rateinv} and
completes the proof of Lemma \ref{lem: rateinv}.




\setcounter{section}{0}
\begin{supplement}\label{sec: proofstechnical}
  \stitle{Supplement of ``Harris recurrent Markov chains and nonlinear monotone cointegrated models''}
  \sdescription{This is the supplementary material associated with the present article.}

  \section{Markov chain theory and notation}\label{sec: extended_markov_theory}
  This section extends Section \ref{sec: markov_theory} of the main paper, presenting a
  more detailed exposition of the Markov chain theory required in the proofs. For further
  details, we refer the reader to \cite{Nummelin1984,Meyn2009,markovChain2018}.

  Let \(\chain = X_0, X_1, X_2, \ldots\) be a time-homogeneous Markov Chain defined on a
  probability space \(\left(E, \SgE,\P\right)\) where \(\SgE\) is countably generated.
  Let \(P\left(x,A\right)\) denote its transition kernel, i.e. for \(x\in E\) ,
  \(A\in\SgE\) we have
  \[P \left( {x,A} \right) = \P\left( {{X_{i + 1}} \in A\left| {{X_i} = x} \right.} \right),\quad i = 0,1, \ldots \]

  Let \(P^n(x,A)\) denote the \(n\)-step transition probability, i.e.
  \[{P^n}\left( {x,A} \right) = \P\left( {{X_{i + n}} \in A\left| {{X_i} = x} \right.} \right)\quad \forall i.\]

  If \(\lambda\) is a probability measure in \(\left(E, \SgE\right)\) such that
  \(\mathcal{L}\left(X_0\right)=\lambda\), then \(\lambda\) is called the \textit{initial
  measure} of the chain \(\chain\). A homogeneous Markov chain is uniquely identified by
  its kernel and initial measure.

  When the initial measure of the chain is given, we will write \(\P_{\lambda}\) (and
  \(\E_{\lambda}\)) for the probability (and the expectation) conditioned on
  \(\mathcal{L}\left(X_0\right)=\lambda\). When \(\lambda=\delta_{x}\) for some \(x\in
  E\) we will simply write \(\P_x\) and \(\E_{x}\).

  An homogeneous Markov chain is \textit{irreducible} if there exists a \(\sigma\)-finite
  measure \(\phi\) on \(\left(E, \SgE\right)\) such that for all \(x\in E\) and all
  \(A\in\SgE\) with \(\phi(A)>0\) we have \(P^n(x,A)>0\) for some \(n\geq 1\). In this
  case, there exists a maximal irreducibility measure \(\irreducibilityMeasure\) (all
  other irreducibility measures are absolutely continuous with respect to
  \(\irreducibilityMeasure\)). In the following, all Markov chains are supposed to be
  irreducible with maximal irreducibility measure \(\irreducibilityMeasure\).

  For any set \(C\in\SgE\), we will denote by \(\sigma_C\) and \(\tau_C\), respectively,
  the times of first visit and first return of the chain to the set \(C\), i.e.
  \(\tau_C=\operatorname{inf}\left\{ n\geq 1: X_n\in C \right\}\) and
  \(\sigma_C=\operatorname{inf}\left\{ n\geq 0: X_n\in C \right\}\). The subsequent visit
  and return times \(\sigma_C, \tau_C\left(k\right)\), \(k\geq 1\) are defined
  inductively as follows:
  \begin{align}
    \tau_C\left(1\right) = \tau_C\quad     & ,\quad \tau_C\left(k\right) = \operatorname{min}\left\{ n>\tau_C\left( k-1 \right):X_n\in C\right\}\label{eq: ch1: time_return_k}, \\
    \sigma_C\left(1\right) = \sigma_C\quad & ,\quad \sigma_C\left(k\right) = \operatorname{min}\left\{ n>\sigma_C\left( k-1 \right):X_n\in C\right\}.
  \end{align}

  Given that our methods will only deal with the values of \(\chain\) in a fixed set
  \(C\), if \(A\) is a measurable set, we will write \(\I_C\{X_t\in A\}\) instead of
  \(\I\{X_t\in A\cap C\}\) and if \(A=E\), then we will simply write
  \(\I_C\left(X_t\right)\).

  We will use \(T_n\left(C\right)\) to denote the random variable that counts the number
  of times the chain has visited the set \(C\) up to time \(n\), that is \(T_n\left( C
  \right)=\sum_{t=0}^{n}{\I_C\left(X_t\right)}\). Similarly, we will write
  \(T\left(C\right)\) for the total of numbers of visits the chain \(\chain\) to \(C\).
  The set \(C\) is called \textit{recurrent} if \(\E_x T\left(C\right)=+\infty\) for all
  \(x\in C\). The chain \(\chain\) is considered recurrent if every set $A\in\SgE$, such
  that \(\irreducibilityMeasure\left(A\right)>0\), is recurrent.

  Although recurrent chains possess many interesting properties, a stronger type of
  recurrence is required in our analysis. An irreducible Markov chain is \textit{Harris
  recurrent} if for all \(x\in E\) and all \(A\in\SgE\) with
  \(\irreducibilityMeasure(A)>0\) we have
  \[\P\left( {{X_n} \in A\; \text{infinitely often}\left| {{X_0} = x} \right.} \right) = 1.\]

  An irreducible chain possesses an accessible atom, if there is a set \(\atom\in\SgE\)
  such that for all \(x,y\) in \(\atom\): \(P(x,.)=P(y,.)\) and
  \(\irreducibilityMeasure(\atom)>0\). Denote by \(\P_{\atom}\) and
  \(\expectation{\atom}{.}\) the probability and the expectation conditionally to
  \(X_0\in \atom\). If \(\chain\) possesses an accessible atom and is Harris recurrent,
  the probability of returning infinitely often to the atom \(\atom\) is equal to one, no
  matter the starting point, i.e. \(\forall x\in E,{\P}_{x}\left(
  \tau_{\atom}<\infty\right) =1.\) Moreover, it follows from the \textit{strong Markov
  property} that the sample paths may be divided into independent blocks of random length
  corresponding to consecutive visits to \(\atom\):
  \begin{align*}
    \block_0 & =\left( X_0,X_1,\ldots,X_{\tau_{\atom}\left( 1 \right)} \right)                              \\
    \block_1 & =\left( X_{\tau_{\atom}\left( 1 \right)+1},\ldots,X_{\tau_{\atom}\left( 2 \right)} \right)   \\
             & \ldots                                                                                       \\
    \block_n & =\left( X_{\tau_{\atom}\left( n \right)+1},\ldots,X_{\tau_{\atom}\left( n+1 \right)} \right) \\
             & \ldots
  \end{align*}
  taking their values in the torus \(\mathbb{T}=\cup_{n=1}^{\infty}E^{n}\). Notice that the
  distribution of \(\block_0\) depends on the initial measure, therefore it does not have the
  same distribution as \(\block_j\) for \(j\geq 1\).
  The sequence \(\left\{\tau_{\atom}(j)\right\} _{j\geqslant 1}\) defines successive times
  at which the chain forgets its past, called \textit{regeneration times}. Similarly, the
  sequence of i.i.d. blocks \(\{\block_j\}_{j\geq 1}\) are named \textit{regeneration blocks}.
  The random variable \(\numreg=T_n\left(\atom\right)-1\) counts the number of i.i.d. blocks up
  to time \(n\). This term is called \textit{number of regenerations up to time \(n\)}.

  Notice that for any function defined on \(E\), we can write
  \(\sum_{t=0}^{n}{f\left(X_t\right)}\) as a sum of independent random variables as
  follows:
  \begin{equation}\label{eq:suplement:split_sum_equation}
    \sum_{t=0}^{n}{f\left(X_t\right)}=f\left(\block_0\right)+\sum_{j=1}^{T\left( n \right)}{f\left( \block_j \right)}+{f\left( \block_{(n)} \right)},
  \end{equation}
  where, \(f\left( \block_0 \right)=\sum_{t=0}^{\tau_{\alpha}}{f\left( X_t \right)}\),
  \(f\left( \block_j \right)=\sum_{t=\tau_{\atom}\left( j \right)+1}^{\tau_{\alpha}\left( j+1 \right)}{f\left( X_t \right)}\)
  for \(j=1,\ldots, \numreg\) and ${f\left( \block_{(n)} \right)}$ denotes the last incomplete block, i.e.
  \({f\left( \block_{(n)} \right)}=\sum_{t=\tau_{\atom}\left( T\left( n \right)+1 \right)+1}^n{f\left(X_t\right)}\).

  When an accessible atom exists, the \textit{stochastic stability} properties of
  \(\chain\) amount to properties concerning the speed of return time to the atom only.
  For instance, the measure \(\pi_{\atom}\) given by:
  \begin{equation}\label{eq:suplement:invariantMeasure}
    \pi_{\atom}\left(B\right) = \E_{\atom} \left( \sum_{n=1}^{\tau_{\atom}} \I{\left\{X_{i}\in B\right\}}  \right),
    \quad \forall B\in\SgE
  \end{equation}
  is invariant, i.e.
  \[\pi_{\atom} \left( B \right) = \int {P\left( {x,B} \right)d\pi_{\atom} \left( x \right)}.\]

  Denote by \(\SgE^+\) the class of nonnegative measurable functions with positive
  \(\irreducibilityMeasure\) support. A function \(s\in\SgE^+\) is called \textit{small}
  if there exists an integer \(m_0\geq 1\) and a measure \(\nu\in\measures\) such that
  \begin{equation}\label{eq:suplement:minorization_general}
    {P ^{m_0}}\left( {x,A} \right) \geqslant s\left( x \right)\nu \left( A \right)\quad\forall x\in E, A\in\SgE.
  \end{equation}
  When a chain possesses a small function \(s\), we say that it satisfies the
  \textit{minorization inequality} \(M\left(m_0,s,\nu\right)\). As pointed out in
  \cite{Nummelin1984}, there is no loss of generality in assuming that
  \(0\leq s\left(x\right)\leq 1\) and  \(\int_{E}{s(x)d\nu(x)}>0\).

  A set \(A\in\SgE\) is said to be \textit{small} if the function \(\I_A\) is small.
  Similarly, a measure \(\nu\) is \textit{small} if there exist \(m_0\), and \(s\) that
  satisfy \eqref{eq:suplement:minorization_general}. By Theorem 2.1 in
  \cite{Nummelin1984}, every irreducible Markov chain possesses a small function and
  Proposition 2.6 of the same book shows that every measurable set \(A\) with
  \(\irreducibilityMeasure\left(A\right)>0\) contains a small set. In practice, finding
  such a set consists in most cases in exhibiting an accessible set, for which the
  probability that the chain returns to it in $m$ steps is uniformly bounded below.
  Moreover, under quite wide conditions a compact set will be small, see
  \cite{Feigin1985}.

  If \(\chain\) does not possess an atom but is Harris recurrent (and therefore satisfies
  a minorization inequality \(M\left(m_0,s,\nu\right)\)), a \textit{splitting technique},
  introduced in \cite{Nummelin1978, Nummelin1984}, allows us to extend in some sense the
  probabilistic structure of \(\chain\) in order to artificially construct an atom. The
  general idea behind this construction is to expand the sample space so as to define a
  sequence $(Y_{n})_{n\in\mathbb{N}} $ of Bernoulli r.v.'s and a bivariate chain
  \(\splitchain=\left\{ {\left( {{X_n},{Y_n}} \right)} \right\}_{n=0}^{+\infty}\), named
  \textit{split chain}, such that the set \(\check{\atom}=\left(E,1\right)\) is an atom
  of this chain. A detailed description of this construction can be found in
  \cite{Nummelin1984}.

  The whole point of this construction consists in the fact that $\splitchain$ inherits
  all the communication and stochastic stability properties from $\chain$
  (irreducibility, Harris recurrence,...). In particular, the marginal distribution of
  the first coordinate process of \(\splitchain\) and the distribution of the original
  \(\chain\) are identical. Hence, the splitting method enables us to establish all the
  results known for atomic chains to general Harris chains, for example, the existence of
  an invariant measure which is unique up to multiplicative constant (see Proposition
  10.4.2 in \cite{Meyn2009}).

  The invariant measure is finite if and only if
  \(\E_{\check{\atom}}{\tau_{\check{\atom}}}<+\infty\), in this case we say the chain is
  \textit{positive recurrent}, otherwise, we say the chain is \textit{null recurrent}. A
  null recurrent chain is called \(\beta\)-null recurrent (c.f. Definition 3.2 in
  \cite{Tjostheim-2001}) if there exists a small nonnegative function \(h\), a
  probability measure \(\lambda\), a constant \(\beta \in \left( {0,1} \right)\) and a
  slowly varying function \(L_h\) such that
  \[{\E_\lambda }\left( {\sum\limits_{t = 0}^n {h\left( {{X_t}} \right)} } \right) \sim \frac{1}{{\Gamma \left( {1 + \beta } \right)}}{n^\beta }{L_h}\left( n \right)\quad\textnormal{as}\;n\to\infty.\]

  As argued in \cite{Tjostheim-2001}, is not a too severe restriction to assume
  \(m_0=1\). Therefore, throughout this paper we assume that \(\chain\) satisfies the
  minorization inequality \(M(1,s,\nu)\), i.e, there exist a measurable function \(s\)
  and a probability measure \(\nu\) such that \(0\leq s\left(x\right)\leq 1\),
  \(\int_{E}{s(x)d\nu(x)}>0\) and
  \begin{equation}\label{eq:suplement:minorization}
    {P}\left( {x,A} \right) \geqslant s\left( x \right)\nu \left( A \right).
  \end{equation}

  A measurable and positive function \(L\), defined in \([a,+\infty)\) for some \(a\geq
  0\), is called \textit{slowly varying at \(+\infty\)} if it satisfies \(\lim_{x\to
  +\infty}{\frac{L\left( xt \right)}{L\left( x \right)}}=1\) for all \(t\geq a\). See
  \cite{Bingham1987} for a detailed compendium of these types of functions.

  It was shown in Theorem 3.1 of \cite{Tjostheim-2001} that if the chain satisfies the
  minorization condition (\ref{eq: minorization}), then it is \(\beta\)-null recurrent if
  and only if
  \begin{equation}\label{def:suplement:betaNullRecurrent}
    \P\left( {{\tau_{\check{\atom}} } > n} \right) \sim \frac{1}{{{\Gamma\left(1-\beta\right) n^\beta }L\left( n \right)}},
  \end{equation}
  where \(L\) is a slowly varying function.

  The following theorem is a compendium of the main properties of Harris's recurrent
  Markov chains that will be used throughout the proofs. Among other things, it shows
  that the asymptotic behavior of \(\numreg\) is similar to the function
  \(u\left(n\right)\) defined as
  \begin{equation}\label{eq:suplement:definition_un}
    u\left(n\right)=\begin{cases}
      n,                       & \text{if }\chain \text{ is positive recurrent}          \\
      n^\beta L\left(n\right), & \text{if }\chain \text{ is }\beta\text{-null recurrent}
    \end{cases}.
  \end{equation}

  The Mittag-Leffler distribution with index \(\beta\) is a non-negative continuous
  distribution, whose moments are given by
  \begin{equation*}
    \E \left(M^m_\beta\left(1\right)\right)=\frac{m!}{\Gamma\left(1+m\beta\right)}\;\;m\geq 0.
  \end{equation*}
  See (3.39) in \cite{Tjostheim-2001} for more details.

  \begin{theorem}\label{th:suplement:t_n_general_results} Suppose \(\chain\) is a Harris recurrent, irreducible Markov chain, with initial measure \(\lambda\), that
    satisfies the minorization condition \eqref{eq:suplement:minorization}. Let \(\numreg\)
    be the number of complete regenerations until time \(n\) of the split chain \(\splitchain\)
    , let \(C\in\SgE\) be a small set and \(\pi\) be an invariant measure for \(\chain\).
    Then,
    \begin{enumerate}
      \item \(0<\pi\left(C\right)<+\infty\).
      \item For any function \(f\), defined on \(E\), the decomposition
            \eqref{eq:suplement:split_sum_equation} holds. Moreover, there is a constant
            \(K_\pi\), that only depends on \(\pi\), such that if \(f\in\measurables\),
            then \(\E_{\lambda}{f\left(\block_1\right)}=K_\pi\int_E fd\pi\).
      \item \(\frac{\numreg}{\Tn}\) converges almost surely to a positive constant.
      \item \(\frac{\numreg}{u\left(n\right)}\) converges almost surely to a positive
            constant if \(\chain\) is positive recurrent
            and converges in distribution to a Mittag-Leffler
            random variable with index \(\beta\) if \(\chain\) is \(\beta\)-null recurrent.
    \end{enumerate}
  \end{theorem}
  \begin{remark}
    The Mittag-Leffler  distribution with index \(\beta\)
    is a non-negative continuous distribution, whose moments
    are given by
    \begin{equation*}
      \E \left(M^m_\beta\left(1\right)\right)=\frac{m!}{\Gamma\left(1+m\beta\right)}\;\;m\geq 0.
    \end{equation*}
    See (3.39) in \cite{Tjostheim-2001} for more details.
  \end{remark}
  \begin{remark}
    Part 1 of Theorem \ref{th:suplement:t_n_general_results} is Proposition 5.6.ii of \cite{Nummelin1984}, part 2 is
    equation (3.23) of \cite{Tjostheim-2001} and part 3 is an application of the Ratio Limit Theorem (Theorem 17.2.1 of \cite{Meyn2009}). For the positive recurrent case, part 4 also follows by
    the aforementioned Ratio Limit Theorem while the claim for the null recurrent case appears
    as Theorem 3.2 in \cite{Tjostheim-2001}.
  \end{remark}

  \section{Technical proofs for Section \ref{sec: localized_markov_chains}}\label{sec: prooflocalized_markov_chains}

  \par {\bf Proof of Lemma \ref{lem: cvFn}.}
  Equation \eqref{eq: cvFnloc} follows from Corollary 2 in \cite{athreya2016general} and part 2 of
  Theorem \ref{th:suplement:t_n_general_results}.

  Now, we turn to the proof of \eqref{eq: cvFninvloc}. To do this, we adapt some of the
  ideas presented in the proof of Lemma 21.2 in \cite{van2000asymptotic}.

  Let $V$ be a normal random variable independent of the $X_i$'s, and $\Phi$ its
  distribution function. it follows from \eqref{eq: cvFnloc} that conditionally on the
  \(X_t\)'s, $F_n(V)$ converges almost surely to $\distfunc(V)$. Thus, denoting by
  \(\P_X\) the conditional probability given the \(X_t\)'s, it follows from \eqref{eq:
  invFn} that $\Phi(F_n^{-1}(u))=\P_X(F_n(V)<u)$ converges almost surely to
  $\P_X(\distfunc(V)<u)=\Phi(\distfunc^{-1}(u))$ at every $u$ at which the limit function
  is continuous . Since $\distfunc$ is strictly increasing in \(C\), one can find
  $\epsilon>0$ such that $\distfunc^{-1}$ is continuous on
  $[\distfunc(x_0)-\epsilon,\distfunc(x_0)+\epsilon]$, so the above limit function is
  continuous at every $u\in [\distfunc(x_0)-\epsilon,\distfunc(x_0)+\epsilon]$. By
  continuity of $\Phi^{-1}$ on $(0,1)$, $F_n^{-1}(u)$ converges almost surely to
  $\distfunc^{-1}(u)$ for every such $u$. By monotonicity, the convergence is uniform,
  hence
  \begin{equation*}
    \sup_{\vert p-\distfunc(x_0)\vert\leq\epsilon}\vert F_n^{-1}(p)-\distfunc^{-1}(p)\vert=o(1)\quad a.s.
  \end{equation*}
  as $n\to\infty$.
  \par{\bf Proof of Lemma \ref{lem: covering+VC}.}
  This proof is an adaptation to the localized case of the proof of Lemma 2 in \cite{bertail-portier:2019}. Let $f^\prime_C\in \mathcal F^\prime_C$, i.e., there exists $f\in \mathcal F$ such that $f_C^\prime(B) = \int f(y) \, M_C(B,\mathrm{d}
    y)$. By Cauchy–Schwarz inequality,
  \[{\left( {\int f (y){\mkern 1mu} {M_C}(B,dy)} \right)^2} \leqslant {\ell_C}\left( B \right)\left( {\int {{f^2}{M_C}(B,dy)} } \right),\]
  then
  \[\E_{Q'}(f'^2_C) \leq \mathbb E_{Q'} \left( \ell_C(B)  \left(\int  f (y)^2 \, M_C(B,d y)\right) \right) = \E_{Q_C}(f^2)E_{Q'}(\ell_C^2),\]
  where the last equality follows from \eqref{eq:integral_q}. Applying this to the
  function
  \begin{align*}
    f_C'(B) - f_k'(B) = \int  (f (y)-f_k(y)) \, M_C(B,d y),
  \end{align*}
  when each $f_k$ is the center of an $\epsilon$-cover of the space $\mathcal F$ and $\|f-f_k\|_{L_2(Q_C)}\leq \epsilon$ gives the first assertion of the lemma.
  To obtain the second assertion, note that $U_C' = U\ell_C$ is an envelope for $\mathcal F_C^\prime$. In addition, we have that
  \begin{align*}
    \|U_C'\|_{L_2(Q')} = U \|\ell_C\|_{L_2(Q')} .
  \end{align*}
  From this, we derive that, for every $0<\epsilon<1$,
  \begin{align*}
    \mathcal N (\epsilon \|U_C'\| _{L_2(Q')},\, \mathcal U_C',\, L_2(Q') ) =  \mathcal N (\epsilon U  \|\ell_C\|_{L_2(Q')},\, \mathcal U',\, L_2(Q') ).
  \end{align*}
  Then using the first assertion of the lemma, we obtain for every $0<\epsilon<1$,
  \begin{align*}
    \mathcal N (\epsilon \|U_C'\|_{L_2(Q')},\, \mathcal F_C',\, L_2(Q') )\leq \mathcal N (\epsilon U,\, \mathcal F,\, L_2(Q_C) ),
  \end{align*}
  which implies the second assertion of the Lemmaz whenever the class $\mathcal F$ is VC with envelope $U$.
  \par{\bf Proof of Lemma \ref{lem: cvFnInvRate}.}
  Let \(B\in \hat{E}\) and \(g:E\times\R\to\R_+\). For each \(y\in\R\)
  we define \(g_y\left(x\right)=g\left(x,y\right)\), then, using the notation
  of section \ref{sec: prelimrate} we will have
  \(\hat{g}_y\left(B\right) = \sum_{x\in B\cap C}g\left(x,y\right)\). Finally,
  for any function \(h:\R\to\R\), we define
  \[\tilde{g}^h_y\left(B\right)= (\widehat{g_y-h(y)})(B)= \sum_{x\in B\cap C}\left(g\left(x,y\right) - h\left(y\right)\right)=\hat{g}_y\left(B\right)-\ell_C\left(B\right)h\left(y\right).\]

  Let \(g\left(x,y\right)=\I\left\{x\leq y\right\}\), and \(h=\distfunc\) as defined in
  \eqref{eq: F}. Then, \(\hat{g}_y\left(B\right)=\sum_{x\in B}\I_C{\left\{ x\leq y
  \right\}} \) and
  \[\tilde{g}_y^{\distfunc}\left(B\right)=\sum_{x\in B\cap C}\left(\I\left\{x\leq y\right\} - \distfunc\left(y\right)\right)=\hat{g}_y\left(B\right)-\ell_C\left(B\right)\distfunc\left(y\right).\]

  From now on, we'll remove the superindex from \(\tilde{g}_y^{\distfunc}\) to ease the
  notation.

  By the definition of \(F_n\) and \(\distfunc\) (\eqref{eq: Fn} and \eqref{eq: F}), we
  have that
  \begin{align*}
    F_n(y)-\distfunc(y) & =\frac 1{\Tn} \sum_{i=1}^{\Tn}{\left(\I{\{ X_{\sigma_C\left(i\right)}\leq y \}}-\distfunc(y)\right)}\nonumber                                                               \\
                        & = \frac{1}{\Tn}\sum_{i=0}^{n}{\left(\I_C{\{ X_{t}\leq y \}}-\I_C{\{X_i\} }\distfunc(y)\right)}\nonumber                                                                     \\
                        & = \frac{1}{\Tn}\left( \tilde{g}_y\left( \block_0 \right) + \sum_{i=1}^{T\left( n \right)}{\tilde{g}_y\left(\block_i\right)} + \tilde{g}_y\left(\block_{(n)}\right) \right),
  \end{align*}
  therefore,
  \begin{equation*}
    \sqrt{\Tn} \Bigl( F_n(y)-\distfunc(y) \Bigr) =\frac{\tilde{g}_y\left( \block_0 \right)}{\sqrt{\Tn}} + \frac{\sum_{i=1}^{T\left( n \right)}{\tilde{g}_y\left(\block_i\right)}}{\sqrt{\Tn}} + \frac{\tilde{g}_y\left(\block_{(n)}\right)}{\sqrt{\Tn}}.
  \end{equation*}

  Notice that \(\left|\tilde{g}_y\left( \block_0\right)\right| \leq
  2\ell_C(\block_0)<+\infty\) and \(\Tn\to +\infty\) almost surely, therefore, the first
  term in the last equation converges almost surely to 0 uniformly in \(y\). For the last
  term, we have that
  \[
    \frac{\left|\tilde{g}_y\left(\block_{(n)}\right)\right|}{\sqrt{\Tn}}\le \frac{2\ell_C(\block_{T\left( n \right)})}{\sqrt{\Tn}}=2\sqrt{\frac{T\left( n \right)}{\Tn}}\frac{\ell_C(\block_{T\left( n \right)})}{\sqrt{T\left( n \right)}},
  \]
  by \AssumpSizeBlock, the expectation of \(\ell^2_C(\block_1)\) is finite, then, Lemma 1
  in \cite{athreya2016general} shows that \(\frac{\ell^2_C(\block_n)}{n}\to 0\) a.s.
  which implies that \(\frac{\ell_C(\block_{n})}{\sqrt{n}}\) also converges to 0 a.s.
  Since \(T\left( n \right)\to +\infty\) a.s., by Theorem 6.8.1 in \cite{Gut2013} we have
  \(\frac{\ell_C(\block_{T\left( n \right)})}{\sqrt{T\left( n \right)}}\to 0\) almost
  surely. Joining this with the almost sure convergence of \(\frac{T(n)}{\Tn}\) to a
  positive constant (see Theorem \ref{th:suplement:t_n_general_results}) we obtain that
  \(\frac{\left|\tilde{g}_y\left(\block_{(n)}\right)\right|}{\sqrt{\Tn}}\) converges
  almost surely to 0 uniformly in \(y\). Therefore,

  \begin{equation}\label{eq: decompositionProcess}
    \sqrt{\Tn} \Bigl( F_n(y)-\distfunc(y) \Bigr) =\frac{\sum_{i=1}^{T\left( n \right)}{\tilde{g}_y\left(\block_i\right)}}{\sqrt{T(n)}} + o_P\left( 1 \right).
  \end{equation}

  where we have used that \(\frac{\Tn}{T(n)}\) converges almost surely to a positive
  constant to use \(T(n)\) instead of \(\Tn\).

  Then, \eqref{eq: approxP} will be proved if we show that, for \(\epsilon\) small enough
  \begin{equation}\label{eq: showTightProcess}
    \mathop {\sup }\limits_{|y - x_0| \leq \epsilon } {\frac{\left| \sum_{i=1}^{T\left( n \right)}{\tilde{g}_y\left(\block_i\right)} \right|}{\sqrt{T(n)}}} = {O_p}\left( 1 \right).
  \end{equation}

  Fix \(\eta>0\) arbitrarily. By Lemma \ref{lem: Tn} and Slutsky's theorem, we can find
  positive numbers \(\underline{a}_\eta, \overline{a}_\eta\) and an integer \(N_\eta\)
  such that \(P\left(\mathcal E_n\right)\geq 1-\frac{\eta}{2}\) for all \(n\geq N_\eta\),
  where
  \begin{equation}\label{eq: e_n_set}
    \mathcal E_n= \left\{ \underline{a}_\eta \asymptnumbreg \leq T(n)\leq \overline{a}_\eta \asymptnumbreg\right\}.
  \end{equation}

  Define \(W_n(\epsilon)=\mathop {\sup }\limits_{|y - x_0|\leq \epsilon }{\left|
  \sum_{i=1}^{n}{\tilde{g}_y\left(\block_i\right)} \right|}\) and let \(M_\eta\) be a
  fixed positive number. Then, for all $n\geq N_\eta$
  \begin{align}
    \P\left( \frac{1}{\sqrt{T(n)}}W_{T(n)}>M_\eta \right) & < \P\left( \left\{ \frac{1}{\sqrt{T(n)}}W_{T(n)}>M_\eta \right\}\cap \mathcal E_n\right) + 1-\P\left( \mathcal E_n\right)\nonumber \\
                                                          & <\P\left( \left\{ \frac{1}{\sqrt{T(n)}}W_{T(n)}>M_\eta \right\}\cap \mathcal E_n\right) + \frac{\eta}{2}.\label{eq: prob_bound}
  \end{align}

  On \(\mathcal E_n\), \(\underline{a}_\eta \asymptnumbreg \leq T(n)\leq
  \overline{a}_\eta \asymptnumbreg\), therefore for all \(n\geq N_\eta\)
  \begin{align}
    \P\left( \left\{ \frac{1}{\sqrt{T(n)}}W_{T(n)}>M_\eta \right\}\cap \mathcal E_n\right) & < \P\left( \left\{ \frac{1}{\sqrt{\underline{a}_\eta \asymptnumbreg}} \max_{1\leq k\leq \overline{a}_\eta \asymptnumbreg}{W_{k}}>M_\eta \right\}\cap \mathcal E_n\right),\nonumber \\
                                                                                           & <\P\left( \frac{1}{\sqrt{\underline{a}_\eta \asymptnumbreg}} \max_{1\leq k\leq \overline{a}_\eta \asymptnumbreg}{W_{k}}>M_\eta \right).\label{eq: prob_bound_2}
  \end{align}

  The random variables \(\left\{\tilde{g}_{(\cdot)} \left( \block_k
  \right)\right\}_{k=1}^{\overline{a}_\eta \asymptnumbreg}\) are i.i.d., therefore, by
  Montgomery-Smith's inequality (Lemma 4 in \cite{Adamczak2008}), there exists a
  universal constant \(K\) such that for all \(n\geq N_\eta\),
  \begin{align}
    \P\left( \frac{1}{\sqrt{\underline{a}_\eta \asymptnumbreg}} \max_{1\leq k\leq \overline{a}_\eta \asymptnumbreg}{W_{k}}>M_\eta \right) & < K \P\left( \frac{1}{\sqrt{\underline{a}_\eta \asymptnumbreg}}{W_{\overline{a}_\eta \asymptnumbreg}}>\frac{M_\eta}{K} \right),\nonumber                                                                                                                                   \\
                                                                                                                                          & <K \P\left( \frac{1}{\sqrt{\underline{a}_\eta \asymptnumbreg}}{ \mathop {\sup }\limits_{|y - x_0|\leq \epsilon }{\left| \sum_{i=1}^{ \overline{a}_\eta \asymptnumbreg }{\tilde{g}_y\left(\block_i\right)} \right|} }>\frac{M_\eta}{K} \right)\label{eq: montgomery_smith}.
  \end{align}

  For an arbitrary set \(T\), let \(\ell^{+\infty}(T)\) be the space of all uniformly
  bounded, real functions on \(T\), equipped with the uniform norm. Weak convergence to a
  tight process in this space is characterized by asymptotic tightness plus convergence
  of marginals (see Chapter 1.5 in \cite{vanDerVaart1996}).

  The class of functions \(\mathcal{G}-\distfunc= \left\{ g_y(\cdot) -
  \distfunc\left(y\right)\right\}_{y\in\R}\) is VC with constant envelope \(2\), hence,
  by Lemma \ref{lem: covering+VC}, the class of functions \(\hat{\mathcal{G}-\distfunc}\)
  is also VC and has \(2\ell_C\) as envelope. \(\E \ell_C^2(\block_1)\) is finite (by
  \AssumpSizeBlock), therefore, by Theorem 2.5 in \cite{Kosorok2008},
  \(\hat{\mathcal{G}-\distfunc}\) is Donsker. Then, the process
  \(\frac{1}{\sqrt{\underline{a}_\eta \asymptnumbreg}}{\left| \sum_{i=1}^{
  \overline{a}_\eta \asymptnumbreg }{\tilde{g}_y\left(\block_i\right)} \right|}\)
  converges weakly in \(\ell^\infty\left[ \hat{\mathcal{G}-\distfunc} \right]\) to a
  tight process \(Z\). The map \(y\mapsto \left\| y \right\|_\infty\) from
  \(\ell^\infty\left[ \hat{\mathcal{G}-\distfunc} \right]\) to \(\R\) is continuous with
  respect to the supremum norm (cf. pp 278 of \cite{van2000asymptotic}), therefore,
  \(\frac{1}{\sqrt{\underline{a}_\eta \asymptnumbreg}}{ \mathop {\sup }\limits_{|y -
  x_0|\leq \epsilon }{\left| \sum_{i=1}^{ \overline{a}_\eta \asymptnumbreg
  }{\tilde{g}_y\left(\block_i\right)} \right|} }\) converges in distribution to \({
  \mathop {\sup }\limits_{|y - x_0|\leq \epsilon }}{Z\left( y \right)}\), hence, we can
  find \(V_\eta\) and \(N^\prime_\eta\) such that
  \begin{equation}\label{eq: probSup}
    \P\left( \frac{1}{\sqrt{\underline{a}_\eta \asymptnumbreg}}{ \mathop {\sup }\limits_{|y - x_0|\leq \epsilon }{\left| \sum_{i=1}^{ \overline{a}_\eta \asymptnumbreg }{\tilde{g}_y\left(\block_i\right)} \right|} }>V_\eta \right)<\frac{\eta}{2K},\quad\forall n>N^\prime_\eta.
  \end{equation}
  Choosing \(M_\eta=K V_\eta\) in \ref{eq: probSup} and joining \eqref{eq: montgomery_smith}, \eqref{eq: prob_bound_2} and \eqref{eq: prob_bound}, completes the proof of \eqref{eq: approxP}.

  Now we proceed to prove \eqref{eq: approxinvP}. Let \(\eta\) be fixed, by \eqref{eq:
  approxP} and Lemma \ref{lem: Tn}, we can find \(\epsilon^\prime\), \(M_\eta^\prime\)
  and \(N_\eta^\prime\) such that
  \begin{align}
    \P\left( \sqrt{\Tn}\mathop {\sup }\limits_{|y - x_0| \leq \epsilon^\prime } {\left| {F_n(y) - \distfunc(y)} \right|>M_\eta^\prime}  \right) & <\frac{\eta}{4}\quad\forall n\geq N_\eta^\prime\label{eq: boundedInProb} \\
    \P\left( \mathcal{D}_n \right)                                                                                                              & \geq 1-\frac{\eta}{2}\quad\forall n\geq N_\eta^\prime\label{eq: probsDn}
  \end{align}
  where \(\mathcal D_n= \left\{ \underline{a}_\eta \asymptnumbreg \leq \Tn\leq \overline{a}_\eta \asymptnumbreg\right\}\). Define the sets
  \begin{align*}
    U_n   & = \left\{ \sqrt{\Tn}\mathop {\sup }\limits_{|p - \distfunc(x_0)| \leq \epsilon } {\left| {F_n^{ - 1}(p) - {\distfunc^{ - 1}}(p)} \right|} >M_\eta \right\}, \\
    U_n^1 & = \left\{ \exists p\in \left[ \distfunc(x_0)-\epsilon, \distfunc(x_0)+\epsilon \right]: F_n^{-1}(p)-\distfunc^{-1}(p) >\frac{M_\eta}{\sqrt{\Tn}}\right\},   \\
    U_n^2 & = \left\{ \exists p\in \left[ \distfunc(x_0)-\epsilon, \distfunc(x_0)+\epsilon \right]: \distfunc^{-1}(p)-F_n^{-1}(p) >\frac{M_\eta}{\sqrt{\Tn}}\right\}.
  \end{align*}
  where \(\epsilon\) and \(M_\eta\) are constants that will be specified later.

  On \(U_n^1\cap\mathcal{D}_n\), \(F_n^{-1}\left( p
  \right)>\frac{M_\eta}{\sqrt{\Tn}}+\distfunc^{-1}\left( p
  \right)>\frac{M_\eta}{\sqrt{\overline{a}_\eta \asymptnumbreg}}+\distfunc^{-1}\left( p
  \right)\), hence,
  \begin{align}\label{eq: boundInvFn}
    F_n\left( \frac{M_\eta}{\sqrt{\overline{a}_\eta \asymptnumbreg}}+\distfunc^{-1}\left( p \right) \right) & \leq F_n\left( F_n^{-1}\left( p \right) \right)\leq p+\frac{1}{\Tn}\nonumber \\
                                                                                                            & \leq p+\frac{1}{\underline{a}_\eta \asymptnumbreg}.
  \end{align}
  Assumption {\AssumpDerivee} indicates that \(\distfunc\) has bounded derivative
  in \(C\), take \(K_1\) as the maximum value of this derivative in \(C\), then, the Mean Value Theorem implies that
  \begin{align*}
    p & =\distfunc\left( \distfunc^{-1}\left( p \right) \right)=\distfunc\left( \frac{M_\eta}{\sqrt{\overline{a}_\eta \asymptnumbreg}}+\distfunc^{-1}\left( p \right)  \right)-\frac{\distfunc^\prime\left( \theta_p \right)M_\eta}{\sqrt{\overline{a}_\eta \asymptnumbreg}}\nonumber \\
      & \leq \distfunc\left( \frac{M_\eta}{\sqrt{\overline{a}_\eta \asymptnumbreg}}+\distfunc^{-1}\left( p \right)  \right)-\frac{K_1 M_\eta}{\sqrt{\overline{a}_\eta \asymptnumbreg}}.
  \end{align*}

  After plugging this into \eqref{eq: boundInvFn} we get
  \begin{equation*}
    \distfunc\left( \frac{M_\eta}{\sqrt{\overline{a}_\eta \asymptnumbreg}}+\distfunc^{-1}\left( p \right)  \right)-F_n\left( \frac{M_\eta}{\sqrt{\overline{a}_\eta \asymptnumbreg}}+\distfunc^{-1}\left( p \right) \right)\geq \frac{K_1 M_\eta}{\sqrt{\overline{a}_\eta \asymptnumbreg}}-\frac{1}{\underline{a}_\eta \asymptnumbreg}.
  \end{equation*}
  Because \(\asymptnumbreg\to +\infty\), we can find \(N_1\) such that \(\sqrt{\frac{\overline{a}_\eta}{\asymptnumbreg}}\frac{1}{\underline{a}_\eta K_1}<1\) for all \(n\geq N_1\), taking
  \(M_\eta\) bigger than \(\frac{M_\eta^\prime}{K_1}\sqrt{\frac{\overline{a}_\eta}{\underline{a}_\eta}}+1\)
  and using that \(\Tn\leq \underline{a}_\eta \asymptnumbreg\) on \(\mathcal{D}_n\), we obtain,
  for all $n\geq N_1$
  \begin{equation}\label{eq: boundDiffFcFn}
    \distfunc\left( \frac{M_\eta}{\sqrt{\overline{a}_\eta \asymptnumbreg}}+\distfunc^{-1}\left( p \right)  \right)-F_n\left( \frac{M_\eta}{\sqrt{\overline{a}_\eta \asymptnumbreg}}+\distfunc^{-1}\left( p \right) \right) > \frac{M_\eta^\prime}{\sqrt{\Tn}}.
  \end{equation}
  Let \(N_{2,\eta}\) be such that \(\frac{M_\eta}{\sqrt{\overline{a}_\eta \asymptnumbreg}}<\frac{\epsilon^\prime}{2}\) for \(n\geq N_{2,\eta}\).
  By the continuity of \(\distfunc^{-1}\) in \(\distfunc(x_0)\) there exists
  \(\epsilon>0\) such that \(\left| \distfunc^{-1}\left( p \right) - x_0 \right|\leq\frac{\epsilon^\prime}{2}\) for all \(p\) in \(\left[ \distfunc(x_0)-\epsilon, \distfunc(x_0)+\epsilon \right]\), therefore, the triangular inequality implies that
  \(\frac{M_\eta}{\sqrt{\overline{a}_\eta \asymptnumbreg}}+\distfunc^{-1}\left( p \right)\) lies
  in the interval \(\left[ x_0-\epsilon^\prime,x_0+\epsilon^\prime \right]\)
  for all \(n\geq N_\eta=\text{max}\left( N_1, N_{2,\eta} \right)\). This, alongside \eqref{eq: boundDiffFcFn}, shows that for all \(n\geq N_\eta\)
  \begin{align*}
    U_n^1 \cap\mathcal{D}_n & \subseteq \left\{ \exists y\in \left[ x_0-\epsilon^\prime, x_0+\epsilon^\prime \right]: \distfunc(y)-F_n(y) >\frac{M^\prime_\eta}{\sqrt{\Tn}}\right\} \\
                            & \subseteq \left\{ \sqrt{\Tn}\mathop {\sup }\limits_{|y - x_0| \leq \epsilon^\prime } {\left| {F_n(y) - \distfunc(y)} \right|>M_\eta^\prime} \right\}.
  \end{align*}
  By a similar argument, it can be shown that
  \begin{equation*}
    U_n^2 \cap\mathcal{D}_n\subseteq \left\{ \exists y\in \left[ x_0-\epsilon^\prime, x_0+\epsilon^\prime \right]: F_n(y)-\distfunc(y) >\frac{M^\prime_\eta}{\sqrt{\Tn}}\right\}\quad\forall n\geq N_\eta.
  \end{equation*}

  Using \eqref{eq: boundedInProb} and \(U_n = U_n^1\cup U_n^2\) we obtain that \(\P\left(
  U_n\cap \mathcal{D}_n \right)\leq \frac{\eta}{2}\) for all \(n\geq N_\eta\). Equation
  \eqref{eq: approxinvP} now follows by \eqref{eq: probsDn}.

\end{supplement}


\bibliographystyle{imsart-number} 
\bibliography{bibliography}       

\end{document}